\documentclass[11pt, a4paper]{article}

\usepackage{natbib}
 \bibpunct[, ]{(}{)}{,}{a}{}{,}%
\usepackage{amssymb,amsmath,graphicx}
\usepackage{float}
\usepackage{afterpage}
\usepackage{multirow}
\usepackage{multicol}
\usepackage{tabularx}
\usepackage{hhline}
\usepackage{color}
\usepackage[normalem]{ulem}
\usepackage{dirtytalk}
\usepackage{booktabs}
\usepackage[makeroom]{cancel}
\usepackage{soul}
\usepackage{verbatim}
\usepackage{graphicx}
\usepackage{subfig}
\usepackage{adjustbox}

\usepackage{algorithm}
\usepackage{algpseudocode}
\usepackage{authblk}
\usepackage[a4paper, margin=1in]{geometry}
\usepackage{tabularx}
\usepackage[hyphens]{url}
\usepackage[hidelinks]{hyperref}

\newcommand{\cG}{\mathcal G}
\newcommand{\cA}{\mathcal A}
\newcommand{\cN}{\mathcal N}
\newcommand{\term}{\Theta}
\newcommand{\ctrain}{\Sigma}
\newcommand{\cB}{\mathcal B}
\newcommand{\cBa}{\mathcal B_{a}}
\newcommand{\cBk}{\mathcal B_{k}}

\newcommand{\cK}{\mathcal K}
\newcommand{\cT}{\mathcal T}
\newcommand{\cKr}{\cK_{b}}

\newcommand{\ot}{o_{\sigma}}
\newcommand{\at}{\alpha_{\sigma}}
\newcommand{\dt}{d_{\sigma}}
\newcommand{\bt}{\beta_{\sigma}}

\newcommand{\itermt}{\Theta_{\sigma}^{\textsc{int}}}
\newcommand{\itrain}{i({\sigma})}
\newcommand{\ait}{\alpha_{\itrain}}
\newcommand{\bit}{\beta_{\itrain}}

\newcommand{\cNtrainin}{\cN_{\sigma}^{\textsc{tin}}}
\newcommand{\cNtrainout}{\cN_{\sigma}^{\textsc{tout}}}
\newcommand{\cNtin}{\cN^{\textsc{tin}}}
\newcommand{\cNtout}{\cN^{\textsc{tout}}}
\newcommand{\cNtinterm}{\cN_{\theta}^{\textsc{tin}}}
\newcommand{\cNtoutterm}{\cN_{\theta}^{\textsc{tout}}}
\newcommand{\ti}{t(i)}

\newcommand{\cAtrainm}{\cA_{\sigma}^{\textsc{tm}}}
\newcommand{\cAtrainh}{\cA_{\sigma}^{\textsc{th}}}
\newcommand{\ua}{u_{a}}

\newcommand{\cNbdterm}{\cN_{\theta}^{\textsc{bd}}}
\newcommand{\cNbtterm}{\cN_{\theta}^{\textsc{bt}}}
\newcommand{\cNboterm}{\cN_{\theta}^{\textsc{bo}}}
\newcommand{\cAbb}{\cA^{\textsc{bb}}}
\newcommand{\cAba}{\cA^{\textsc{ba}}}
\newcommand{\cAbd}{\cA^{\textsc{bd}}}
\newcommand{\cAbt}{\cA^{\textsc{bt}}}

\newcommand{\ob}{o_{b}}

\newcommand{\ab}{\alpha_{b}}

\newcommand{\db}{d_{b}}

\newcommand{\bb}{\beta_{b}}

\newcommand{\cBtermiout}{\cB^{+}_{\theta i}}
\newcommand{\cBtermiin}{\cB^{-}_{\theta i}}

\newcommand{\ub}{u_{b}}
\newcommand{\fb}{f_{b}}

\newcommand{\ok}{o_{k}}
\newcommand{\ak}{\alpha_{k}}
\newcommand{\dk}{d_{k}}
\newcommand{\bk}{\beta_{k}}
\newcommand{\tk}{\tau_{k}}

\newcommand{\cNdin}{\cN^{\textsc{din}}}
\newcommand{\cNdout}{\cN^{\textsc{dout}}}
\newcommand{\cNdinterm}{\cN_{\theta}^{\textsc{din}}}
\newcommand{\cNdintermk}{\cN_{\ok}^{\textsc{din}}}
\newcommand{\cNdoutterm}{\cN_{\theta}^{\textsc{dout}}}
\newcommand{\cAdin}{\cA^{\textsc{din}}}

\newcommand{\cAdwt}{\cA^{\textsc{dwt}}}

\newcommand{\cAdout}{\cA^{\textsc{dout}}}

\newcommand{\cNcoterm}{\cN_{\theta}^{\textsc{co}}}
\newcommand{\cNcdterm}{\cN_{\theta}^{\textsc{cd}}}

\newcommand{\cNpoolterm}{\cN^{\textsc{pool}}_\theta}
\newcommand{\cAcb}{\cA^{\textsc{cb}}}
\newcommand{\cAcd}{\cA^{\textsc{cd}}}
\newcommand{\cAcdp}{\cA^{\textsc{up}}} 
\newcommand{\cAbdp}{\cA^{\textsc{ep}}}

\newcommand{\cApvl}{\cA^{\textsc{pl}}}
\newcommand{\cApeb}{\cA^{\textsc{pb}}}

\newcommand{\cApvltermi}{\cA^{\textsc{pl}}_{\theta i}}

\newcommand{\cApebtermi}{\cA^{\textsc{pb}}_{\theta i}}

\newcommand{\cAbdptermi}{\cA^{\textsc{ep}}_{\theta i}}

\newcommand{\cAcdptermi}{\cA^{\textsc{up}}_{\theta i}}

\newcommand{\yb}{y_{b}}

\newcommand{\zrk}{z_{bk}}
\newcommand{\zk}{z_k}
\begin{document}

%\setcounter{page}{1}
%\DOI{10.1287/xxxx.0000.0000}%

%\RUNAUTHOR{Kienzle et al.}

%\RUNTITLE{Intermodal Railroad Blocking and Fleet Management}

% \TITLE{Bundling Information Goods of Decreasing Value}
\title{The Intermodal Railroad Blocking and Railcar Fleet-Management Planning Problem}

\author[1]{Julie Kienzle}
\author[2]{Serge Bisaillon}
\author[3]{Teodor Gabriel Crainic}
\author[4]{Emma Frejinger}
\affil[1]{\small{Department of Computer Science and Operations Research, Universit\'e de Montr\'eal and Interuniversity Research Centre on Enterprise Networks Logistics and Transportation (CIRRELT)}}
\affil[2]{\small{Interuniversity Research Centre on Enterprise Networks Logistics and Transportation (CIRRELT)}}
\affil[3]{\small{Department of Analytics, Operations, and Information Technologies, Universit\'{e} du Qu\'{e}bec \`{a} Montr\'{e}al and Intereruniversity Research Centre on Enterprise Networks Logistics and Transportation (CIRRELT)}}
\affil[4]{\small{Department of Computer Science and Operations Research, Universit\'e de Montr\'eal and Interuniversity Research Centre on Enterprise Networks Logistics and Transportation (CIRRELT)}}

\maketitle

\begin{abstract}
Rail is a cost-effective and relatively low-emission mode for transporting intermodal containers over long distances. This paper addresses tactical planning of intermodal railroad operations by introducing a new problem that simultaneously considers three consolidation processes and the management of a heterogeneous railcar fleet. 
We model the problem with a scheduled service network design with resource management (SSND-RM) formulation, expressed as an 
%path-based 
integer linear program. While such formulations are challenging to solve at scale, we demonstrate that our problem can be tackled with a general-purpose solver when provided with high-quality warm-start solutions. To this end, we design a construction heuristic inspired by a relax-and-fix procedure.
We evaluate the methodology on realistic, large-scale instances from our industrial partner, the Canadian National Railway Company: a North American Class I railroad. The computational experiments show that the proposed approach efficiently solves practically relevant instances, and that solutions to the SSND-RM formulation yield substantially more accurate capacity estimations  compared to those obtained from simpler baseline models.
Managerial insights from our study highlight that ignoring railcar fleet management or container loading constraints can lead to a severe underestimation of required capacity, which may result in costly operational inefficiencies. Furthermore, our results show that the use of multi-platform railcars improves overall capacity utilization and benefits the network, even if they can locally lead to less efficient loading as measured by terminal-level slot utilization performance indicators.
\end{abstract}


\section{Introduction} \label{intro}

Railroads move large quantities of a wide range of valuable commodities over long distances. Rail is considered a cost-effective and relatively low-emission mode for long-distance ground transportation and can therefore play an important role in enhancing the sustainability of freight transportation.

Railroads are a key component of national and global supply chains and \emph{intermodal transportation} that is steadily and strongly growing in volume and worth.
%networks, which has experienced steady traffic growth.  Rail freight transportation can be broadly divided into two main markets: general cargo and intermodal.
Rail intermodal transportation means moving containers loaded in single or double stacks on particular railcars.
Given its economic importance, intermodal transport has become a major component of railroad operations, planning and, management.
%General cargo is transported by cars designed to carry unpackaged goods such as coal, grain, oil, and lumber. The intermodal market consists of moving containers that are placed onto railcars. 
It is well known that, in order to operate efficiently and profitably, railroads face challenging decision-making problems at all levels of planning, strategic, tactical, and operational. 
%Whereas general cargo and intermodal traffic share the same infrastructure, their train services and equipment requirements differ. Therefore, 
Intermodal services require particular attention in this context due to both customer-related factors (e.g., the need to synchronize rail–ship activities in ports and stricter-than-usual time requirements) and operational concerns (e.g., double-stack loaded intermodal railcars generally bypass automated sorting in classification yards). As a result, these services typically call for dedicated planning processes.

We focus on tactical planning of intermodal railroad operations in the North American market. This is a large market with infrastructure that enables double-stacking containers on railcars. While this increases the number of containers that can be transported within a given train length, it also makes the planning problem more complex. Indeed, there are several container and railcar types with specific loading constraints \citep{Mantovani2017}. We introduce the \emph{Intermodal Railroad Blocking and Railcar Fleet-Management} (\emph{IBRM}) problem. Taking a train service schedule as given by the railroad, we aim to build an economically and customer-service efficient plan that simultaneously considers the selection of extra services, the loading of containers on railcars, the blocking of loaded and empty railcars, the selection of blocks, the distribution of demand flows through the service network, and the management of the railcar fleet. The plan is built for a cyclic schedule of given length (e.g., a week), intended to be repeatedly executed over the tactical-planning horizon (e.g., a season).

Particularly challenging issues in addressing the IBRM are (i) simultaneously considering three consolidation processes: containers to railcars, railcars to blocks, and blocks to trains; (ii) differentiating railcar and container types and representing, in a computationally efficient way, the loading of containers onto railcars within a tactical model; and (iii) integrating the management of a limited heterogeneous railcar fleet with the design of the blocking and service network.

While there is an extensive literature on railroad planning, relatively few studies focus on integrated planning and intermodal rail transport.
%Moreover, for tractability, interrelated decision-making problems are often solved in sequence. For example, blocking is addressed first, followed by service selection, with resource management treated as a separate planning problem. 
For example, \cite{Morganti2020} is the study closest to our work, as they focus on intermodal block planning; however, they do not consider railcar fleet management nor extra train service selection.

This paper offers several contributions. \emph{First}, we introduce the IBRM problem, which is new to the literature, and propose an enhanced \emph{Scheduled Service Network Design with Resource Management} (\emph{SSND-RM}) model \citep{crainic+h2021SND_Book}. The SSND-RM model takes the form of an integer linear programming model (ILP) based on a cyclic four-layer space time network representation (namely, the container, car, block and train layers). This approach enables the use of a continuous-time network representation, where the time structure is defined by the arrival and departure times of the train services considered at the terminals on their respective routes. 

\emph{Second}, we propose a construction heuristic inspired by a relax-and-fix procedure. We use the heuristic to compute high-quality warm-start solutions for general-purpose ILP solvers. This simple approach is appealing from a practical point of view, as a general-purpose solver can be used as a black box. Using this approach, we show that we can solve large-scale realistic instances from our industrial partner, the Canadian National Railway Company, in a reasonable time. 

\emph{Third}, by comparing solutions with those obtained by solving simpler baseline models, we show the significance of the problem we introduce. Based on 
%carefully designed 
extensive numerical 
%results 
experiments and 
%the related 
analyses, we offer several managerial insights. In brief, ignoring loading constraints and railcar fleet management can result in severe underestimation of required capacity; we link characteristics of demand (number of containers of different types) to an adequate composition of the railcar fleet; and we highlight the importance of multi-platform cars from a capacity usage point of view.

The structure of the paper unfolds as follows: Section~\ref{BlockCar} describes the problem in detail and Section~\ref{litReview} gives a review of the related literature. We describe the network representation in Section~\ref{network} followed by the SSND-RM formulation in Section~\ref{model}. Section~\ref{solApp} outlines the solution approach, and Section~\ref{application} describes the computational results based on large-scale instances from our industrial partner. Finally, Section~\ref{conclude} concludes the paper. 

\section{Problem Description} \label{BlockCar}

We consider a tactical planning problem facing railroads. The tactical plan is cyclic over a given schedule length $T$, for example, a week, and is repeated over a medium-term planning horizon (e.g., three months). The aim is to guide intermodal rail operations over a physical network. To describe the problem in detail, we start by introducing core concepts -- scheduled train service supply, container demand, and railcars -- followed by a description of loading rules, blocks, and the decision-making problem. We summarize the related notation in Table~\ref{tab:system-notation}.

\paragraph{Train services.} Railroads operate an infrastructure network made up of rail tracks connecting terminals $\theta \in \term$ that are dedicated, totally or partially, to intermodal traffic. We consider a set of scheduled train services $\sigma \in \Sigma$. The services are dedicated to intermodal traffic and each train service $\sigma$ has an origin terminal $\ot$ with departure time $\at$, and a destination terminal $d_{\sigma}$ with arrival time $\bt$. 
It also has a set of intermediate stops $\itermt$, each $i(\sigma) \in \itermt$ with arrival $\bit$ and departure $\ait$ times.
We define a leg as the route between two consecutive stops.
%Each
The service 
%also has maximum 
capacity is expressed in train length (we use feet, shorthand ft) that may vary 
%in our application 
%between different terminals, so-called train 
by leg. We note that intermodal traffic is typically light compared to general cargo, the limiting factor is therefore train length rather than tonnage.  

We consider two types of scheduled train services: First, regular services, $\Sigma_\textsc{init}$, 
%that are operated in each time period and are hence 
%New NOTE Theo: ALL selected services are operated at ALL periods
designed to transport the majority of the recurrent demand. Second, potential extra train services, $\Sigma_\textsc{extra}$, that are to be selected, 
%scheduled, 
%New NOTE Theo: Theu are already scheduled !!!
and operated at fixed cost, $f_{\sigma}, \sigma \in \Sigma_\textsc{extra}$, when needed.
%there is sufficient demand. The latter is expressed as a
An extra service is activated only when loaded at least at
%when minimum percentage of capacity usage 
$U_\sigma \%,~\sigma \in \Sigma_\textsc{extra}$, 
of its capacity on its leg set. 
Extra services can include using capacity on general cargo trains, or dedicated extra intermodal trains.
The scheduled services $\ctrain = \Sigma_\textsc{init} \cup \Sigma_\textsc{extra}$ are given, extra ones being potential only.
%New NOTE Theo: we plan for the horizon !!, nothing is "updated" in the deterministi version, neither regular, nor extras, nor demand or resources.
%taken as given in our problem and they are updated on a less frequent basis than the intermodal plan we focus on. 
The selection of a subset of extra trains is therefore part of the tactical planning problem.

%\textbf{NOTE Theo: I modified the previous paragraph because the set Singma-extra is a set of potential services out of which we must select. I think the modified wording reflects this. Moreover, the explanation on when things are revised is very general and I do not believe is part of the problem setting. If you want to have such an explanation, we should also indicate that the set Sigma-init is determined at a different decision planning process when the entire service network, general and intermodal, is selected.}

\paragraph{Containers.} Intermodal containers come in different sizes and types. For example, the North American market has high- and low-cube containers of various sizes (e.g., 20-, 40-, and 53-ft). Unlike operational load planning problems \citep{Mantovani2017}, only high-level information is available at the tactical planning stage. As we further detail below, crucial in this context are aspects that impact the usage of railcars, and hence the capacity usage of the train. 
Let $k \in \cK$ denote intermodal demand, defined as a number $\upsilon_k$ of containers of a specific type $\tk \in \cT$ arriving at origin terminal $\ok$ at time $\ak$, loaded on railcars and moved to their destination terminal $\dk$ before a certain due-date $\bk$. We consider two container types $\cT = \{\tau_{40}, ~\tau_{53}\}$, as 
two 20-ft containers can be defined as one 40-ft-long box, while containers longer than 40 ft (e.g., 45 ft) may be treated as 53-ft-long ones as they follow similar loading rules.

\paragraph{Railcars.} The North American fleet contains several different railcar types. \cite{Mantovani2017} discuss these in depth. Here, we describe the information that is relevant to tactical planning. The types of railcars $\gamma \in \Gamma$ differ in the number of platforms, 1, 3, or 5, and their respective size (length) $\pi \in \Pi$, and can accommodate single- or double-stacked containers. The most common platform types in the North American market are 40 ft and 53 ft long. Each platform is shaped like a well to lower the center of mass of loaded containers. The space is therefore physically restricted by the platform size. Let $\eta^\gamma_{\pi}$ denote the number of platforms of type $\pi \in \Pi$ on railcar type $\gamma \in \Gamma$ (all the platforms that make up a railcar are of the same type), and let $\Gamma_\pi$ denote the set of railcar types made up of platform type $\pi \in \Pi$. Let $\lambda_\gamma$ be the length of railcar type $\gamma \in \Gamma$. 

For the tactical planning problem, it is important to note that, for a given platform size $\pi$, the ratio $\lambda_\gamma / \eta^\gamma_{\pi}$ decreases with increasing number of platforms. This means that, for example, a five-platform railcar uses less train capacity per platform, than a one-platform railcar. However, railcars with more than one platform have additional loading restrictions. We attend to those next. 

\paragraph{Loading Rules.} Each platform may be double-stacked (two container slots, one at the bottom and one on top) or single-stacked (one slot) and loading rules dictate which configurations of container types are allowed on which railcar type \citep{Mantovani2017}. Here, we focus on rules pertaining to container size, as this is the only information available at the tactical planning level. 

First we consider the bottom slot as the loading rule simply dictates that the platform needs be of sufficient size. That is, independently of the railcar type, a 40-ft container can be loaded in the bottom slot on either a 40-ft or 53-ft platform, whereas a 53-ft container requires a 53-ft platform. 

\begin{figure}[htbp]
	\begin{center}
		\includegraphics[width=\textwidth]{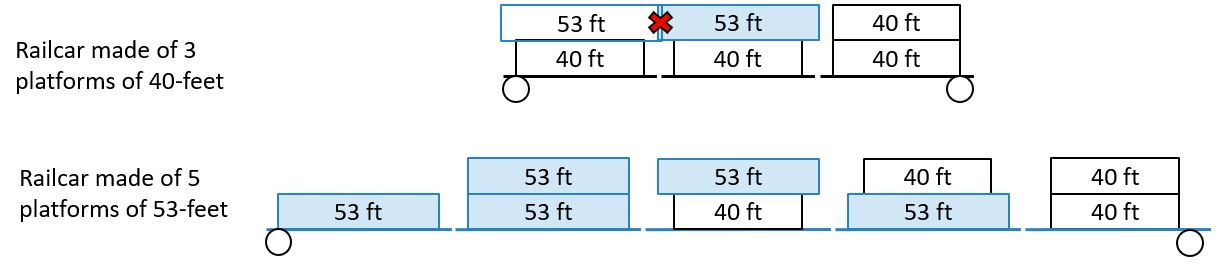}
		\caption{Illustrative examples of loading rules (not respected in the top figure, respected in the bottom figure)}
		\label{fig:railcar} 
	\end{center}
\end{figure}

The top slot size is not physically restricted by a well. So in principle, a 53-ft or 40-ft container could be loaded in the top slot, independently of the container size loaded in the bottom slot. This is true for 40-ft railcars with one platform and all 53-ft railcar types (see the bottom part of Figure~\ref{fig:railcar} for an example of a loading that respects the rules). However, for 40-ft railcars with more than one platform, the limited distance between platforms does not allow to load 53-ft containers in all top slots. Therefore, 40-ft containers can be loaded in any top slot, but 53-ft containers can only be loaded in every second slot (see the top part of Figure~\ref{fig:railcar} for an example of a loading that does not respect the rules). This means that the maximum number of 53-ft containers on railcars with three 40-ft platforms is two, the corresponding figure for five 40-ft platforms is three, etc.

These loading rules imply that several demands with possibly different container types and time characteristics can be loaded together on a same railcar as long as the loading rules are respected. 
	
\paragraph{Blocks.} A block $b \in \cB$ is a consolidation of railcars that move as a single unit between a pair of terminals taking advantage of economies of scale and reducing terminal handling costs. Blocks are assigned to a sequence of train services. Blocking of intermodal traffic differs from its general cargo counterpart. Whereas it is customary to sort %general cargo 
most railcars at so-called classification yards, this is less common for intermodal railcars. For instance, double-stack railcars cannot be sorted in hump yards for physical reasons (it could destabilize the load).
 This means that the intermodal railcars that form a given block all have the same destination, which is the destination of the block.  
 
 As such, a block $b \in \cB$ is defined by an origin terminal $\ob$ where it is built by grouping the classified (sorted) railcars assigned to it, the time of departure at the origin terminal $\ab$, the destination terminal $\db$ where it is dismantled, the time of arrival at destination $\bb$, the sequence of trains moving the block from its origin to its destination, and the terminals where they are being transferred from one train to another with the respective arrival and departure times given by the associated train-service schedules. A block has a maximum capacity $\ub$ that is limited by the train capacity on each train leg. 

\begin{figure}[htbp]
	\begin{center}
		\includegraphics[width=\textwidth]{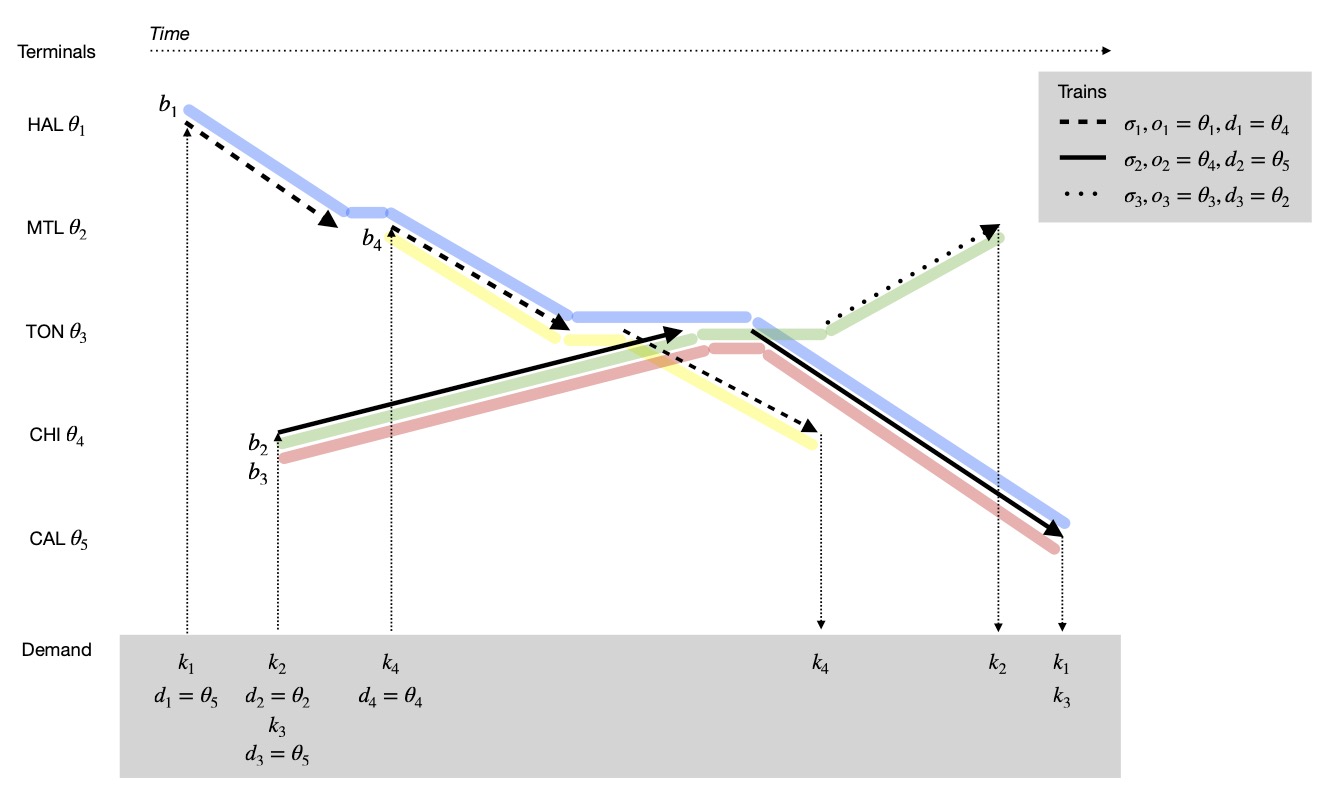}
		\caption{Illustrative example of trains, blocks and demands}
		\label{fig:block} 
	\end{center}
\end{figure}

To illustrate the basic concept, we show a stylized example in Figure~\ref{fig:block}. It depicts three train services $\sigma_1$, $\sigma_2$ and $\sigma_3$ (arrows with dashed, solid, and dotted lines, respectively). There are five terminals represented on the vertical axis ($\theta_1,\ldots,\theta_5$) and time is represented on the horizontal axis. Each train has an origin and a destination terminal 
%(see legend in the figure) a
and are respectively composed of three, two, and one legs. We represent four demands at the bottom part of the figure, each with an origin and a destination. Note that demands $k_2$ and $k_3$ share the same origin, but they have different destinations. For the sake of simplicity, we do not explicitly represent containers and railcars in this figure. The demand is transported on four blocks (see highlighting next to the trains with the block identification at the origin): $b_1$ is transported by trains $\sigma_1$ and $\sigma_2$, $b_2$ by $\sigma_2$ and $\sigma_3$, whereas $b_3$ and $b_4$ are transported by a single train ($\sigma_2$ and $\sigma_1$, respectively).

We note that a given demand $k$ can be split on multiple railcars, and even on multiple blocks. Let $\cBk \subseteq \cB$ denote the set of blocks that can transport demand $k \in \cK$ given its time characteristics, while $\cK_b$ is the set of demands that can be assigned to block $b \in \cB$.
Let $c_{bk}$ be the unit cost of transporting a container of demand $k$ on block $b \in \cBk$, and $c_{b \gamma}$ the unit cost of using railcar type $\gamma$ in block $b \in \cB$. 
The tactical plan selects the subset of blocks in $\cB$ that should be built.
Fixed costs $f_b,~b\in \cB,$ reflect the time and 
%capacity usage 
resources required to build a block at its origin terminal
and to transfer it between services, when relevant. The transportation of demand that cannot be satisfied by the service network defined by the tactical plan is outsourced at a relatively high cost $c_k,~k\in \cK$.
 
\paragraph{Railcar fleet management.} We consider a heterogeneous fleet of railcars where we denote by $H^\gamma$ the total number of railcars of type $\gamma \in \Gamma$. Because there is a limited railcar supply in the network, it is essential to efficiently manage these resources to bring the right type of railcars to the right location in the network at the right time. This is a challenging aspect of the problem as demand may be unbalanced between the origin-destination pairs, and the right mix of railcar types depends on the specific mixes of container types (see discussion on loading rules). 
Some railcars might, thus, have to travel empty so as to be distributed where they are needed (e.g., at terminals where a railcar shortage is occurring). When the fleet composition is modified due to, e.g., seasonal changes or new railcar acquisitions, railroad companies have to decide at which terminals the railcars should be first stationed. In this case, in addition to providing a railcar circulation (i.e., determining the flow of railcars within the network to support the planned services), the railcar fleet management plan should also give the initial allocation of railcars to the different terminals in the network. Let $c_{\gamma \theta}$ denote the unit cost of initially allocating railcar type $\gamma \in \Gamma$ to terminal $\theta \in \Theta$.

\paragraph{Decision-making problem.} Given a train schedule, a set of demands, and a heterogeneous fleet of  railcars, the goal of the IBRM problem is to determine simultaneously a scheduled blocking plan that includes  container-to-railcar (loading problem) and railcar-to-block (blocking problem) assignments, a railcar circulation to support the selected blocks with, when wanted, an initial railcar allocation to terminals (railcar fleet management problem), the demand itineraries, including, when needed, the volumes assigned to general cargo trains or trucks, and finally, when desired, a selection of extra trains to add to the schedule to satisfy demand; all of this while minimizing the total operational costs. As a benefit, managerial strategic insights, such as determining an effective railcar fleet composition, can also be derived from the provided tactical plan. 

\begin{table}[htbp]
    \centering
    \begin{tabular}{l l}
    \hline
    \multicolumn{2}{l}{\textbf{Infrastructure}} \\ 
    $\term$  & Set of terminals, $\theta \in \term$  \\
    $t^{\textsc{trans}}$ & Block transfer time at any terminal \\ 
    \hline
    \multicolumn{2}{l}{\textbf{Train services}  } \\
    $\ctrain$ & Set of scheduled train services, $\sigma \in \ctrain$ \\
    $\Sigma_\textsc{init}$ & Set of regular scheduled train services \\
    $\Sigma_\textsc{extra}$ & Set extra (non regular) scheduled train services \\
    $\ot$ & Origin terminal of $\sigma \in \ctrain$ \\
     $\at$ & Departure time of $\sigma \in \ctrain$ \\
    $\dt$  & Destination terminal of $\sigma \in \ctrain$\\
    $\bt$ & Arrival time of $\sigma \in \ctrain$ \\
    $\itermt$ & Set of intermediate stops of $\sigma \in \ctrain$ \\
    $\bit$ & Arrival time at intermediate stop (node) $i(\sigma)\in \itermt,~\sigma \in \ctrain$ \\ 
    $\ait$ & Departure time at intermediate stop (node) $i(\sigma)\in \itermt,~\sigma \in \ctrain$ \\
    $f_{\sigma}$ & Fixed cost of using $\sigma \in \Sigma_\textsc{extra}$ \\
    $U_\sigma$ & Capacity usage threshold for $\sigma \in \Sigma_\textsc{extra}$  \\
    \hline 
\multicolumn{2}{l}{\textbf{Demand}} \\
    $\cK$ & Set of demands, $k \in \cK$ \\   
    $\ok$ & Origin terminal of $k \in \cK$\\
    $\ak$ & Arrival time at origin of $k \in \cK$\\
    $\dk$ & Destination terminal of $k \in \cK$ \\
    $\bk$ & Due date at destination of $k \in \cK$\\
    $\upsilon_k$ & Number of containers of $k \in \cK$ \\
    $\cT$ & Set of container types, $\tau \in \cT$; $\cT=\{\tau_{40},\tau_{53}\}$ in current setting \\
    $\tk$ & Container type of $k\in \cK$ \\    
    $c_k$ & Unit outsourcing cost $k\in \cK$ \\
    \hline
\multicolumn{2}{l}{\textbf{Railcars}} \\
    $\Gamma$ & Set of railcar types, $\gamma \in \Gamma$ \\
    $\lambda_\gamma$ & Length of railcar type $\gamma \in \Gamma$ \\
    $\Pi$ & Set of platform sizes, $\pi \in \Pi$ $\Pi = \{\pi_{40},\pi_{53}\}$ in current setting \\
    $\Gamma_{\pi}$ & Set of railcar types made up of platform type $\pi \in \Pi$ \\
    $\eta^\gamma_{\pi}$ & Number of platforms of type $\pi$ for railcar $\gamma \in \Gamma$ \\
    $H^\gamma$ & Number of railcars of type $\gamma \in \Gamma$ available in the system\\
    $c_{\gamma \theta}$ & Unit cost for intially allocating railcar type $\gamma \in \Gamma$ to terminal $\theta \in \Theta$  \\
     \hline
\multicolumn{2}{l}{\textbf{Blocks }} \\
    $\cB$ & Set of blocks, $b \in \cB$ \\
    $\cB_k$ & Set of blocks which may transport demand $k\in \cK$, $\cB_k \subseteq \cB$ \\
    $\ob$  & Origin terminal of $b \in \cB$ \\
    $\ab$ & Departure time of $b \in \cB$ \\
    $\db$ & Destination terminal of $b \in \cB$ \\
    $\bb$ & Arrival time of $b \in \cB$ \\
    $\fb$ & Cost of building $b \in \cB$ \\ 
    $\ub$ & Capacity of $b \in \cB$ \\ 
    $\cK_b$ & Set of demands that can be assigned to $b \in \cB$ \\
    $c_{bk}$ & Unit cost of moving a container of demand $k$ on $b \in \cB$ \\
    $c_{b \gamma}$ & Unit cost of moving a railcar of type $\gamma$ on $b \in \cB$\\
   \hline
\end{tabular}
\caption{Summary of notation related to core concepts and problem parameters}
    \label{tab:system-notation}
\end{table}

\section{Literature Review}
\label{litReview}

%\ef{Note not to forget: add literature on heuristics + general purpose solvers... For example, \cite{Fischetti2017} show that adhoc preprocessing on top of a general-purpose solver can be used for real-time train rescheduling. }

There is a rich and long history of successful Operations Research developments and contributions targeting railroad planning.
A comprehensive review of this literature goes far beyond the scope of this paper.
The interested reader may turn to a series of surveys that synthesize this story and contributions, including \citet{assad80b,dejax+cr87,crainic88revrail,crainic+laporte97,Cordeau1998,crainic2000ejor,newman+ny2002,crainic03hall,ahuja+cs2005tutorial,crainic+k07,bektas+c07,crainic09Intermodal,yaghini+a2012ProcediaSNDrev,piu+s2014LocoSurvey} and \citet{Chouman2020}.
The Service Network Design methodology and contributions are synthesized and reviewed in \citet{crainic+h2021SND_Book,crainic2024SNDbasic_ModelBook,crainic2024SNDadvanced_ModelBook} and \citet{crainic+r2024_50Years}.
Several observations and trends stand out when examining this literature and allow one to situate our contribution within the field.

First, few studies target intermodal rail transportation planning. We discuss them in the last part of the section, following a rapid overview of developments relevant to the problem we address.

Second, one observes that most early contributions focus on single problems or combinations of a limited number of them, blocking, service selection, and resource management in particular.
Blocking \citep[see, e.g.,][]{bodin+al80,newton+bv97block,ahuja+jl07,jha+as08,yaghini+a2012ProcediaSNDrev,yaghini+ma2021Networks_Blocking} is assumed to be addressed first \citep[in contrast,][assumes a given service schedule]{Morganti2020}, followed by service selection and makeup \citep[see, e.g.,][]{morlok+peterson70,assad80a,nozick+m97,newman+y2000TS_Schedule,yaghini+a2012ProcediaSNDrev}.
In most cases, network design type formulations \citep{magnanti+w84,crainic+gg2021NetDes_Book} are proposed to select blocks or services, from given sets, respectively.

The management of resources is then addressed as a different planning activity, focusing on distribution and routing given a set of services. To address the challenge of unbalanced trade, resources that complete their assignments and are needed in terminals different from their current locations are repositioned to undertake the next round of activities. 
The activity is often characterized as ``operational'', even though resources support services and, accordingly, the planning horizon of the previous two activities is shared.
In general, network flow optimization methods are proposed in this context, evolving from transportation/assignment models \citep[see, e.g.,][]{bomberault+w66,white68,white+b69,white72} to multicommodity integer flow time-space formulations integrating various practical rules and limitations \citep[see, e.g.,][]{florian+all76,Joborn2004,ahuja+all05,vaidyanathan+ahuja+all08,vaidyanathan+ahuja+o08,balakrishnan+ks2016TS_Crews,belgacem+all2016TS_Engine2Plan,piu+all2015Loco,ortiz+cf2021TS_Loco,miranda+cf2020Loco}.

Third, a trend can be observed in railroad tactical planning toward comprehensive models that integrate the main system and operational components, e.g., selection, scheduling and makeup of services, railcar classification and block design, resource assignment and management, and demand routing through the selected service network.
Service Network Design (SND) appears as the methodology of choice to address these problems \citep{Chouman2020}.
Our work belongs to this railroad planning and OR methodological development area.

\cite{crainic+fr1984TS_rail} propose what is probably the first railroad tactical planning SND model
\citep[][generalizes it for consolidation-based multicommodity multimode freight transportation systems]{crainic+r1986TRB_SND}.
The model integrates service selection, service frequency optimization, car classification and blocking, train makeup, and freight routing.
The nonlinear path-based formulation also accounts for congestion and delays in terminals and on the network's tracks, as well as for the (re-)positioning of empty railcars through one or several OD demand matrices (generated through demand-distribution models from the surplus and penury levels at terminals derived from the loaded demand).

\citet{Hagani1989} presents a SSND model combining train routing, scheduling, and make-up, as well as empty car distribution, on a space-time network.
A heuristic is used to address a somewhat simplified version of the model and illustrate the interest of integrated planning.
The model proposed by \citet{Keaton1989,Keaton1992} aims to determine the pairs of terminals to connect by direct services, whether to offer more than one train a day, as well as the routing of freight and the blocking of railcars.
\citet{Gorman1998} starts from the previous model to design, using a tabu-enhanced genetic search metaheuristic, a scheduled operating plan that follows as much as possible the particular operation rules of a major North-American railroad.
All these contributions model blocking through classification costs, rather than explicit blocking decision variables. 

\citet{zhu+cg2014or} propose what appears to be the first comprehensive SSND railroad planning model, integrating the selection of scheduled services, selection of blocks, service makeup, and railcar classification and routing. 
The model is built on a cyclic, multi-layer \citep{crainic+gk2022MultiLayer,crainic2024MLayerBG_Book} time-space network corresponding to the discretized schedule length. The layers corresponds to the railcar terminal activities (demand entry and exit, and railcar classification and blocking), block terminal activities, and train service operations.
The authors propose a  matheuristic combining dynamic block generation and slope scaling, Tabu Search and ellipsoidal search methods. The authors do not address resource management issues.

\citet{pedersen+c07,pedersen+cm2009TS} focus on integrating the management of one type of resource, namely, locomotives or intermodal shuttles (i.e., locomotive and fixed set of container-carrying railcars), into tactical-planning SSND formulations.
The authors formalize the concept of design-balanced SSND, forcing equal numbers of selected services arriving at or leaving from any terminal, at all relevant time moments.
%These constraints directly regulate the flow of resources when a single unit of resource is assigned to a service, while they are written in terms of number of resource units otherwise.
\citet{pedersen+c07,vu+ct13JOH,chouman+c15TS} present metaheuristics to address the resulting SSND.

\citet{andersen+cc09EJOR} enlarge the scope of the previous problem setting and devise SSND with resource management models to address the coordination and synchronization of services and fleets of several railroads and navigation companies at particular junction points (train services operating in several countries need locomotive changes at some borders due to differences in electric tensions).
\citet{andersen+cc2009TRC} then observe that resources that support services or that are repositioned operate according to cycle schedules, returning to their home bases after certain periods of time. 
The authors compare SSND arc and cycle-based formulations and show that the latter provide more modeling flexibility, in terms of cycle duration for example, and computational efficiency, provided that cycles are generated a priori.
\citet{andersen+gcc2011TS} propose a branch-and-price with column generation algorithm to address the latter issue.
\cite{crainic+htv2014TS,crainic+htv2018EJTL} extend this work in the general SSND context for consolidation-based freight carriers (with one unit of resource per service only) in two ways: 1) a larger range of resource management concerns, e.g., multiple resource types, acquisition and allocation/reallocation of resources at terminals, outsourcing services, and cycle duration rules; 2) a matheuristic that incorporates slope scaling and network flow-based cycle generation heuristics.

To conclude this brief review, we return to the limited number of contributions aimed at planning intermodal railroad transport.
\citet{nozick+m97} minimizes the total moving cost of flatcars and trailers (1 trailer loaded on a flatcar) given a set of services. 
\citet{newman+y2000TS_Schedule} proposes a day-of-week uncapacitated (in the number of trains one can make up in a yard and operate on a line) frequency SND model, to determine whether homogeneous container OD demand should be moved by a direct or indirect (through a unique main yard) service. 
\citet{muller+ee2021CIE_IntermodeStochSSND} propose a SSND with stochastic demand for the case of a homogeneous set of services (called ``vehicles'' in the paper) managed through design-balancing constraints.
A very simple loading rule (one container) characterize these three contributions, which do not consider the selection and formation of services and blocks.

\citet{Morganti2020} start to address the loading issue when formulating their SSND model for the block planning problem for intermodal rail transport with double-stack railcars. Operational loading requires detailed information about containers, such as their individual weights \citep{MantEtAl2016}, and is hence too detailed for tactical planning. 
%The set of rules that govern the operational loading of multiple types of containers on many different intermodal railcars is extremely broad and precise \citep{MantEtAl2016}, much too detailed and complex for tactical planning.
\citet{Morganti2020} therefore propose a first approximation of those rules, based on the length of the main container and railcar platform types, together with constraints enforcing the length of selected blocks and of the services to which they are assigned.
The space-time network is discretized according to the schedules of the given set of services, their respective times of arrival and departure at terminals yielding the time moments of the network. We adopt the same modeling of time in our work.
Unlike us, \citet{Morganti2020} do not consider service selection or resource management.

In conclusion, we note that the current state-of-the-art in railroad tactical planning does not address the main challenges of intermodal systems.
The methodology we propose contributes to filling this gap, integrating within a comprehensive SSND model service and block selection, realistic tactical-level container-to-railcal loading constraints, and the management of multiple sets of railcars that may be combined, loaded or empty, when building the blocks.

\section{Network Modeling} \label{network}

The multi-layer network structure is presented in Section \ref{model:network}, introducing the time representation and the train-service layer. The other layers and the inter-layer arcs are then described by following a demand itinerary, from the initial processing at the origin terminal (Section~\ref{embark}), to moving through the network (Section~\ref{moving}), and the final handling at destination (Section~\ref{disembark}). Recall that the tactical plan defines the demand itineraries, the blocks to build, the extra services to run, and the management of the railcar fleet. Table~\ref{tab:network-notation} summarizes the time-space network notation.

\subsection{General Structure and Train Layer}  \label{model:network}

We propose a four-layer \citep{crainic+gk2022MultiLayer,crainic2024MLayerBG_Book} cyclic time-space network $\cG = (\cN, \cA)$ of length $T$ to represent the system dynamics, the activities, and the decisions of the IBRM problem. Nodes $i \in \cN$ represent events at specific moments in time $\ti$ and terminals $\theta (i)$. Arcs stand for activities taking place in time and space between these nodes.

\begin{figure}[htbp]
	\begin{center}
		\includegraphics[width =\textwidth]{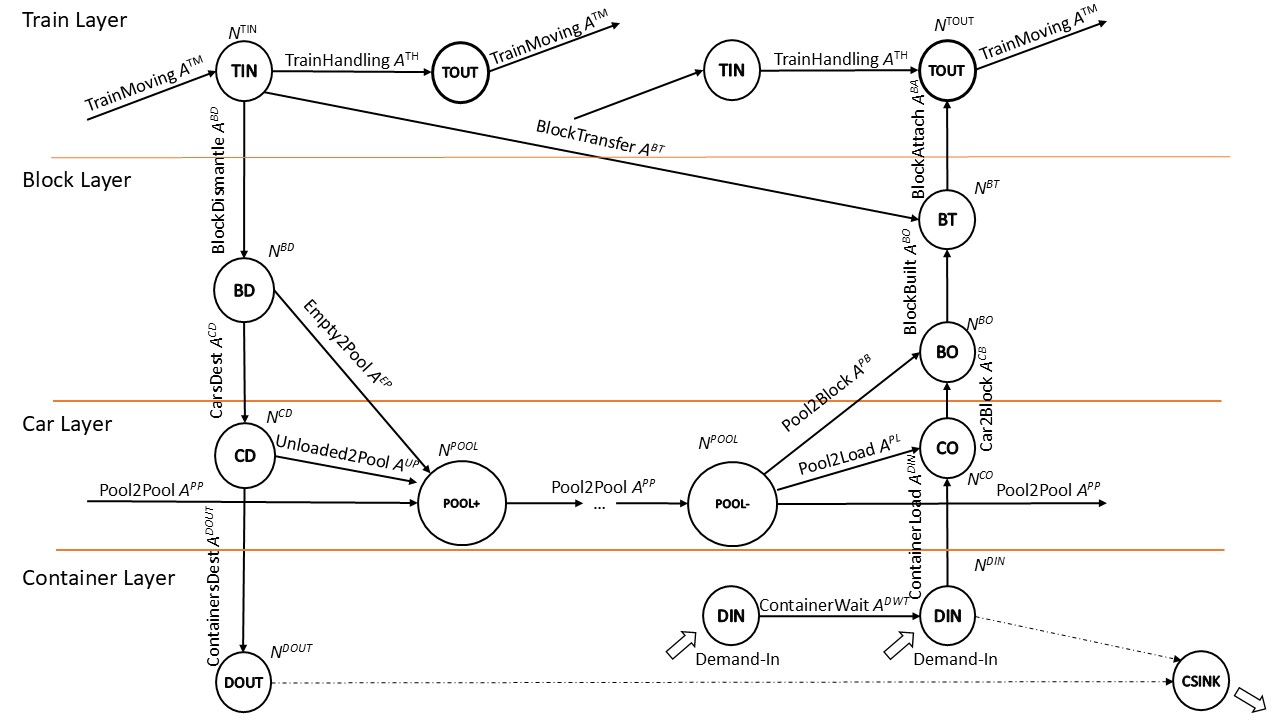}
		\caption{Four-layer (train, block, car, container) time-space network illustration}
		\label{fig:network} 
	\end{center}
\end{figure}

The layers correspond to the main components of the system, train services, blocks, railcars, and containers.
\emph{Intra-layer} nodes and arcs model the activities corresponding to each layer at appropriate moments in time, while \emph{inter-layer} arcs capture the interactions among the components, e.g., loading/unloading of containers onto/from railcars, consolidating the latter into blocks, attaching/detaching the latter to/from train services, and dismantling blocks at destination. Due to the cyclic nature of the tactical plan, activities initiated before the end of the schedule may end after $T$, that is, during the next application of the plan. This is modeled by having the corresponding arcs wraparound, times computed $modulo(T)$ \citep{Chouman2020,crainic+h2021SND_Book}. The network is illustrated in Figure \ref{fig:network} for a terminal $\theta \in \term$ and a representative time interval featuring a train arrival and a train departure. 
We refer to this figure throughout the section.

The \emph{container layer} represents the arrival of each demand $k \in \cK$ at its origin yard, the possible waiting before being loaded onto railcars, the actual loading operation, the unloading at destination, and the subsequent exit from the system. The \emph{car layer} connects the container-based demand to the system and governs how containers are loaded on railcars, that are grouped into blocks and then into trains.
It is also in the \emph{car layer} that the fleets of the various types of railcars are managed through the assignment to container-loading activities, the empty movements to balance the needs within the network, the railcar-pool inventory counts at terminals, as well as the initial fleet size and allocation at terminals. Blocks are handled on the \emph{block layer}: building at origin by consolidating empty and loaded railcars, attachment to a train, transfer between two train services, and detachment from the last train at destination for dismantling. Detailed descriptions of the network components and modeling of these three layers are given in the next subsections.

The \emph{train layer} plays a special role in the proposed formulation.
It obviously represents the movements between terminals and the time spent there to pick up or drop off the blocks of regular and extra train services.
It also defines the time instances of the time-space network and, thus, the pacing of activities and decisions. 
Recall that, most time-space networks proposed in the service network design literature, including models targeting rail planning \citep[see, e.g.,][]{Chouman2020,crainic+h2021SND_Book}, generate time-space networks by 1) discretizing the schedule length according to a given granularity and 2) duplicating each node of the physical network at each of the resulting periods. We take a different approach and define a \emph{continuous-time} network based on the schedules of the given train-service set $\ctrain$.

The arrival and departure times of each train service to/from each terminal on its route yield corresponding arrival (\textsc{tin}) and departure (\textsc{tout}) nodes, defining the time structure of the entire network (Figure \ref{fig:network}). Let $\cNtoutterm$ and $\cNtinterm$ be the sets of \textsc{tout} and \textsc{tin} nodes at terminal $\theta \in \term$, sorted in increasing time order, with $\cNtout = \cup_{\theta \in \term} \cNtoutterm$ and $\cNtin = \cup_{\theta \in \term} \cNtinterm$. Then, any node $i \in \cN$, representing an event on any layer, has its time stamp $t(i)$ equal to either a train arrival or departure moment at terminal $\theta \in \term$, i.e., $i \in \cNtinterm$ or $i \in \cNtoutterm$, respectively. In other words, container, railcar, and block activities are synchronized with train departures and arrivals, the associated waiting being represented on the inter-layer arcs.

Notice that, more than one service could arrive at or leave a terminal at any given time moment, generating several simultaneous nodes. Such an event is rare, however. Hence, to simplify the presentation, and without loss of generality, we assume a single service arrival or departure at any time moment.

\begin{table}[htbp]
    \centering
    \begin{tabularx}{\textwidth}{l X}
    \hline
    $\cG = (\cN, \cA)$ & Time-space network, sets of nodes $\cN$ and arcs $\cA$ \\
    $\theta(i) , \ti$ & Terminal and time moment of node $i \in \cN$ \\
    $T$ & Schedule length \\
    $\ua$ & Capacity of arc $a \in \cA$ \\
    \hline
\multicolumn{2}{l}{\textbf{Train Layer}} \\
     $\cNtin$ & Train-service arrival nodes \& times at terminals (\textsc{tin}), = $\cup_{\theta \in \term} \cNtinterm$  \\
     $\cNtout$ & Train-service departure nodes \& times at terminals (\textsc{tout}),  = $\cup_{\theta \in \term} \cNtoutterm$  \\
     $\cAtrainm$ & \emph{TrainMoving} (\textsc{tout}, \textsc{tin}) moving arcs of service $\sigma \in \Sigma$ \\
     $\mathcal{A}_{\sigma}^{\textsc{utm}}$ & Leg set of extra train service $\sigma \in {\Sigma}_{\textsc{extra}}$ for which a minimum load $U_{\sigma}$ \\ 
     & is required in order for the service to be considered for selection \\
     $\cAtrainh$ & \emph{TrainHandling} (\textsc{tin}, \textsc{tout}) handling arcs at terminals serviced by $\sigma \in \Sigma$ \\ 
     $\cBa$ & Set of blocks on train-service arc $a \in \cAtrainm \cup \cAtrainh$ \\
    \hline
\multicolumn{2}{l}{\textbf{Block Layer}} \\
    $\cNboterm$ & Block forming nodes (\textsc{bo}) \\
    $\cNbdterm$ & Block-at-destination nodes (\textsc{bd}) \\
    $\cNbtterm$ & Block transfer and assembly nodes (\textsc{bt}) \\
    $\cB_\theta^+$ & Blocks formed at
    terminal $\theta \in \term$ \\
    $\cBtermiout$ & Blocks formed at
    terminal $\theta \in \term, t(i) = \ab , i \in \cN^{\textsc{pool}-}_\theta$ \\
    $\cBtermiin$ & Blocks dismanteled at
    terminal $\theta \in \term, t(i) = \bb , i \in \cN^{\textsc{pool}+}_\theta$ \\
    $\cAbb$ & \emph{BlockBuild} arcs moving new blocks to the assembly node (\textsc{bt}) \\
    $\cAbt$ & \emph{BlockTransfer} inter-layer arcs, bring transferred blocks to the assembly node \\ 
    $\cAba$ & \emph{BlockAttach} inter-layer arcs, move collected blocks to the outgoing train \\
    $\cAbd$ & \emph{BlockDismantle} inter-layer arcs, take blocks off trains at destination \\
    \hline
\multicolumn{2}{l}{\textbf{Car Layer}} \\
    $\cN^{\textsc{pool}}$ & Railcar-inventory nodes (\textsc{pool}), $ = \cup_{\theta \in \Theta} \cN^{\textsc{pool}}_{\theta}$ \\
    $\cN^{\textsc{pool}}_{\theta}$ & Terminal $\theta \in \Theta$ railcar-inventory nodes,\\& $=\cN^{\textsc{pool}+}_{\theta} \cup \cN^{\textsc{pool}-}_{\theta}=\cup_{\gamma \in \Gamma} \{ \cN^{\textsc{pool}+}_{\gamma \theta} \cup \cN^{\textsc{pool}-}_{\gamma \theta} \}$ \\
     $\cN^{\textsc{pool}+}_{\gamma \theta}$ & Terminal $\theta \in \Theta$ inventory nodes linked to the arrival of trains\\ 
     & with possibly increasing railcar-type $\gamma \in \Gamma$ inventory\\
     $\cN^{\textsc{pool}-}_{\gamma \theta}$ & Terminal $\theta \in \Theta$ inventory nodes linked to the departure of trains\\ 
     & with possibly decreasing railcar-type $\gamma \in \Gamma$ inventory\\
     $\cNcoterm$ & Container-to-railcar loading nodes (\textsc{co}), $\theta \in \Theta$ \\
     $\cNcdterm$ & Container-from-railcar unloading nodes (\textsc{cd}), $\theta \in \Theta$ \\
     $\cA^{\textsc{pp}}$ & \emph{Pool2Pool} railcar inventory carrying arcs \\
     $\cApvl$ & \emph{Pool2Load} empty-railcar-to-container assignment \& loading arcs \\ 
     $\cApeb$ & \emph{Pool2Block} empty-railcar inter-layer blocking arcs \\ 
     $\cAcb$ & \emph{Cars2Blocks} loaded-railcar inter-layer blocking arcs \\ 
     $\cAbdp$ & \emph{Empty2Pool} empty-railcar at destination inter-layer disassemble arcs \\
     $\cAcd$ & \emph{CarsDest} loaded-railcar at destination inter-layer disassemble arcs \\ 
     $\cAcdp$ & \emph{Unloaded2Pool} arcs moving unloaded (empty) railcars to their inventory pool \\
     \hline
\multicolumn{2}{l}{\textbf{Container Layer}} \\
     $\cNdin$ & Demand-arrival \textsc{din} nodes, $= \cup_{\theta \in \term} \cNdinterm$ \\
     $\cNdintermk \subseteq \cNdinterm$ & Possible \textsc{din} nodes for demand $k \in \cK$ \\
     $\cNdout$ & Demand-at-destination \textsc{dout} nodes, $= \cup_{\theta \in \term} \cNdoutterm$\\
     $\cAdin = \cup_{\theta \in \term} \cAdin_\theta$  & \emph{ContainerLoad} container-to-car origin inter-layer arcs \\
     $\cAdwt$  & \emph{ContainerWait} arcs \\
     $\cAdout$ & \emph{ContainerDest} car-to-container destination inter-layer arcs \\ 
     \hline
    \end{tabularx}
    \caption{Time-space network notation}
    \label{tab:network-notation}
\end{table}

\subsection{Embarking on a Train}  \label{embark}

We describe in the following the ``embarking" side of the container, car, and block layers.
We discuss the initial allocation of the railcar fleets to  terminals in the car layer subsection.

\subsubsection{The Container Layer}
models demand entering (and existing) the system with the set of nodes $\cNdin = \cup_{\theta \in \term} \cNdinterm$. A demand $k \in \cK$ leaves its origin terminal $\ok$ at the earliest when the first train, following its own arrival time at the terminal, departs. As such, a node is generated in the container layer every time a train leaves (instead of every time a demand enters). 
%only (right part of Figure~\ref{fig:demand-node}), rather than every time a demand arrives (left part of Figure \ref{fig:demand-node}).
%the time of each node $i \in \cNdinterm$ (\textsc{din} in Figure~\ref{fig:network}) corresponding to the departure time of a train from terminal $\theta \in \term$.
Accordingly, a demand $k \in \cK$ arriving to $o_k$ at time $\ak$ connects to nodes $i\in \cNdintermk$ such that $ \ti \geq \ak$.
Recall that demands can be divided and delivered using different blocks and trains. This is modeled with a set of \emph{ContainerWait} arcs ${\cAdwt} = \cup_{\theta \in \term}\{ a = (i, j), j=i+1, i \in \cNdinterm  ,\ \theta \in \term \}$ where $j = i+1$ represents the next-in-time DIN node in terminal $\theta$.

%Containers not selected at a \textsc{din} node thus wait until the next loading possibility, i.e., the next-in-time \textsc{din} node. These arcs are collected in the \emph{ContainerWait} arc set ${\cAdwt} = \cup_{\theta \in \term}\{ a = (i  , i+1 \in \cNdinterm)  ,\ \theta \in \term \}$.

The operation of loading demand onto railcars at the origin terminal is represented by inter-layer \emph{ContainerLoad} arcs $\cAdin = \cup_{\theta \in \term} \cAdin_\theta = \{ a = (i,j)|t(i)=t(j), i \in \cNdinterm , j  \in \cNcoterm ,\theta \in \term \}$ linking \textsc{din} nodes in the container layer to container-railcar loading nodes \textsc{co} in the car layer. The container-to-railcar assignments and loading are performed according to the rules introduced in Section~\ref{BlockCar}.

\subsubsection{The Car Layer} Four main activities are modeled in this layer: container loading (and unloading), railcar-to-block assignment, railcar fleet management, and the initial dimensioning of the fleet and its assignment to terminals.

We explicitly represent the railcar circulation within the cyclic scheduled service network to control the inventories and flows of each railcar type. For this purpose, we
introduce \textsc{pool} nodes to model inventory changes to each railcar type at each terminal and relevant time moment.
Railcar inventories at terminals \emph{decrease} when railcars are assigned loaded or empty to a block. We define the associated set of inventory nodes linked to departing trains $\cN^{\textsc{pool}-}_{\theta} = \cup_{\gamma \in \Gamma}\cN^{\textsc{pool}-}_{\gamma \theta} := \{ i ~|~ t(i) = \at, \alpha_{j(\sigma)}, j(\sigma) = 1, \ldots, |\itermt|, \sigma \in \Sigma \}$, 
%for all railcar types 
$\gamma \in \Gamma$,
%and terminals 
$\theta \in \Theta$. 
We model the assignment of loaded and empty railcars to blocks with intra-layer arcs \emph{Pool2Load} ${\cApvl} := \cup_{\theta \in \term}\{ a = (i ,j)|t(i)=t(j), i \in \cN^{\textsc{pool}-}_{\theta} , j \in \cNcoterm ,\ \theta \in \term\}$ 
and inter-layer arcs
%and the assignment of empty railcars to blocks by inter-layer 
\emph{Pool2Block} ${\cApeb} := \cup_{\theta \in \term}\{ a = (i , j)|t(i)=t(j), i \in \cN^{\textsc{pool}-}_{\theta} , j \in \cNboterm ,\theta \in \term\}$, respectively, connecting the pool to the block-forming \textsc{bo} node in the block layer. Let $\cB_\theta^+ = \cup_{i \in \cN^{\textsc{pool}-}_{\theta}} \cBtermiout = \{b \in \cB ~|~ \ob = \theta (i) , \alpha_b = t(i) , i \in \cN^{\textsc{pool}-}_\theta , \theta \in \term \}$ stand for the blocks, containing these empty and loaded railcars, formed at their origin terminal $\theta$ at the corresponding train-departure time.

Analogously, inventories \emph{increase} when blocks are dismantled (see Section~\ref{disembark}), with associated set of inventory nodes, linked to the arrival of a train,  $\cN^{\textsc{pool}+}_{\theta} = \cup_{\gamma \in \Gamma}\cN^{\textsc{pool}+}_{\gamma \theta} = \{ i ~|~ t(i) = \bt, \beta_{j(\sigma)}, j(\sigma) = 1, \ldots, |\itermt|, \sigma \in \Sigma \}$, $\gamma \in \Gamma , \theta \in \Theta$. 
%(We attend to the related activities in Section~\ref{disembark}.)
% and
% $\cN^{\textsc{pool}-}_{\gamma \theta} = \{ i ~|~ t(i) = \at, \alpha_{j(\sigma)}, j(\sigma) = 1, \ldots, |\itermt| \}$, 
% the sets of pool nodes when the inventory of railcar type $\gamma \in \Gamma $ possibly increases or decreases at terminal $\theta \in \Theta$, after the arrival of an incoming train service or when providing for a departing train service,  respectively.
We denote the set of all \textsc{pool} nodes $\cN^{\textsc{pool}}_{\theta} := \cup_{\gamma \in \Gamma} \{ \cN^{\textsc{pool}+}_{\gamma \theta} \cup \cN^{\textsc{pool}-}_{\gamma \theta} \}$ for each terminal $\theta \in \Theta$, and $\cN^{\textsc{pool}} := \cup_{\theta \in \Theta} \cN^{\textsc{pool}}_{\theta}$.

\emph{Pool2Pool} arcs ${\cA^{\textsc{pp}}} := \cup_{\theta \in \term}\{a = (i,j), j = i+1, i,j \in \cN^{\textsc{pool}}_{\theta} \}$ connect successive inventory nodes at each terminal, and carry the flows of available empty railcars.
To ensure the cyclic nature of the service network and schedule, each arc $a \in \cA^{\textsc{pp}}$, starting from the node $i \in \cN^{\textsc{pool}}_{\theta}, \theta \in \term$, with the closest time $t(i)$ to the end of the schedule $T$, wraps around to the first node $i \in \cN^{\textsc{pool}}_{\theta}$ following the starting of the schedule, that is, to the first train-service arrival or departure moment to/from the associated terminal $\theta \in \term$. %Figure \ref{fig:spool-ssink} illustrates this case.

%Empty railcars are taken out of inventory to be loaded and then assigned to a block. \emph{Pool2Load} arcs ${\cApvl} = \cup_{\theta \in \term}\{ a = (i ,j), i \in \cN^{\textsc{pool}}_{\theta} , j \in \cNcoterm ,\ \theta \in \term\}$ represent the former activity, bringing empty railcars from the pool of the terminal to the associated \textsc{co} loading node. 

Inter-layer \emph{Car2Block} arcs ${\cAcb} = \cup_{\theta \in \term}\{ a = (i , j)|t(i)=t(j), i \in \cNcoterm , j \in \cNboterm ,\theta \in \term \}$ model the car-to-block assignment decisions, linking the corresponding nodes in the car (\textsc{co}) and block (\textsc{bo}) layers.

\subsubsection{The Block Layer}

Blocks are built by consolidating empty and loaded railcars brought to a \emph{block-forming} node \textsc{bo} by the Pool2Block ($\cApeb$) and Car2Block ($\cAcb$) arcs, respectively.
The just-built blocks are to be attached to the train leaving at the \textsc{tout} time moment of the blocking operation, together with blocks delivered earlier to the terminal to be transferred to the same train as part of their journeys (see Section \ref{moving}). 
Let \textsc{bt} identify a \emph{block transfer and assembly} node.
The\emph{BlockBuild} arcs ${\cAbb} = \cup_{\theta \in \term}\{ a = (i , j)|t(i)=t(j) , i \in \cNboterm , j \in \cNbtterm) ,\theta \in \term \}$ then bring the newly built blocks to the assembly node, where they are joined by the transferred blocks brought by the inter-layer arcs $\cAbt$ (see Section \ref{moving}).
Inter-layer \emph{BlockAttach} arcs ${\cAba} = \cup_{\theta \in \term}\{ a = (i , j)|t(i)=t(j) , i \in  \cNbtterm , j \in \cNtoutterm) ,\theta \in \term \}$ then move the collected blocks from the assembly node to the \textsc{tout} node of the outgoing train service to which they are assigned.

%Demand flows are now ready to start their journeys towards their destinations.

\subsection{Moving Between Terminals -- The Train Layer} \label{moving}
Two sets of arcs define the train-service operations.
The set of \emph{TrainMoving} arcs $\cAtrainm = \{ a = (i , j) , i \in  \cNtrainout , j \in \cNtrainin) ,\sigma \in \ctrain \}$ link consecutive departure (\textsc{tout}) and arrival (\textsc{tin}) nodes of train services  moving between terminals.
Set $\mathcal{A}_{\sigma}^{\textsc{utm}} \subseteq \cAtrainm , \sigma \in {\Sigma}_{\textsc{extra}}$, contains the legs of an extra train service for which a minimum load is required in order for the service to be considered for selection.
The \emph{TrainHandling} arcs, collected in the set $\cAtrainh = \{ a = (i , j ) , i \in \cNtrainin, j \in \cNtrainout) ,\sigma \in \ctrain \}$, connect consecutive arrival (\textsc{tin}) and departure (\textsc{tout}) nodes of train services, representing the time spent in the terminal to perform block pick-up and drop-off activities (and other maintenance work).

Train services are made up of, and move, blocks.
Let $\cBa , a \in \cAtrainh$, be the set of the blocks that stay on the train service $\sigma \in \ctrain$ while terminal activities take place on its Train-Handling arc $a$. Similarly, let $\cBa , a \in \cAtrainm$, be the set of blocks making up train service $\sigma$ on its Train-Moving arc $a$. The load hauled on each of these arcs $a \in \cA^{\textsc{tm}} = \cup_{\sigma \in \Sigma} \cAtrainm$ is limited by capacity $\ua$. Note that no capacity is enforced on \emph{TrainHandling} arcs as modeling terminal operations is not part of our problem. 

A block is moved by a sequence of train services, its journey corresponding to a path of Train-Moving and Train-Handling arcs, \emph{transfers} taking place when the entire block changes trains at intermediary terminals.  
Inter-layer \emph{BlockTransfer} arcs, gathered in set ${\cAbt} = \cup_{\theta \in \term}\{ a = (i , j) , i \in  \cNtinterm , j \in \cNbtterm) ,\theta \in \term \}$, model this activity linking the train and block layers.
Let $t^{\textsc{trans}}$ be the time required for executing a block transfer. 
\textsc{tin} nodes in the train layer are thus linked to all \textsc{bt} nodes of the corresponding terminal such that 
$t^{\textsc{trans}} \leq t(i) - t(j) , i \in \cNbtterm , j \in \cNtinterm, \theta \in \Theta$.

\subsection{Disembarking From a Train } \label{disembark}

Blocks, and the railcars and containers making them up, are taken off the last train service at the destination terminal and are dismantled.
Railcars carrying containers are then at destination and, once their containers are unloaded, they join the empty railcars which were in the block (if any) in the appropriate inventory pool.
The unloaded containers are delivered to their consignees.
We follow these activities, completing the time-space network definition, in the reverse order of the embarking activities (see Section \ref{embark}).

\subsubsection{The Block Layer}

The blocks arriving at destination are taken off the train and dismantled.
Inter-layer \emph{BlockDismantle} arcs ${\cAbd} = \cup_{\theta \in \term}\{ a = (i , j)|t(i)=t(j) , j \in  \cNtinterm , j \in \cNbdterm) ,\theta \in \term \}$ model these activities, moving the blocks from the train \textsc{tin} nodes to the corresponding block-at-destination \textsc{bd} node at the destination terminal.

Let $\cBtermiin = \{b \in \cB ~|~ \db = \theta (i) , \beta_b = t(i) , i \in \cN^{\textsc{pool}+}_\theta , \theta \in \Theta\}$ stand for the blocks dropped at their destination terminal $\theta$ at the corresponding train-arrival time.
Their dismantlement yields empty railcars, increasing the inventory of the terminal as described next.

\subsubsection{The Car Layer}

Dismantled blocks yield empty and loaded railcars, of various types, which increase their availability at the terminal.
Inter-layer \emph{Empty2Pool} arcs in ${\cAbdp} = \cup_{\theta \in \term}\{ a = (i , j)|t(i)=t(j) , i \in \cNbdterm , j \in \cN^{\textsc{pool}+}_\theta),\theta \in \term \}$ model the former activity, moving empty railcars from block-at-destination \textsc{bd} nodes, in the block layer, directly to the corresponding $\textsc{pool}$ nodes in the car layer.

Loaded railcars must first be unloaded before being added to the inventories corresponding to their types held in the terminal pool.
\emph{CarsDest} inter-layer arcs ${\cAcd} = \cup_{\theta \in \term}\{ a = (i , j)|t(i)=t(j) , i \in  \cNbdterm , j \in \cNcdterm,\theta \in \term \}$ bring the loaded railcars from \textsc{bd} nodes to the railcar-unloading \textsc{cd} nodes in the car layer.
The \emph{Unloaded2Pool} arcs ${\cAcdp} = \cup_{\theta \in \term}\{ a = (i , j)|t(i)=t(j) , i \in \cNcdterm , j \in \cN^{\textsc{pool}+}_\theta,\ \theta \in \term \}$ then bring the now-empty railcars to the pool of the terminal.

\subsubsection{The Container Layer}

Demand flows exit the system through the container layer. The set of  \emph{ContainerDest} inter-layer arcs ${\cAdout} = \cup_{\theta \in \term}\{(i , j)|t(i)=t(j) , i  \in \cNcdterm , j \in \cNdoutterm ,\theta \in \term \}$ groups the arcs bringing unloaded containers from the nodes \textsc{cd} in the car layer, to the $i \in \cNdout$ destination nodes (\textsc{dout}) in the container layer. 
Being in the yard of the destination terminal marks the end of a demand's freight journey.
The final delivery to customers is then modeled by the arc linking the \textsc{dout} node to the container-sink node \textsc{csink} collecting all demand flows sent to the terminal. 

The time-space network is completed by the artificial arcs capturing the demand volumes not assigned to services in $\ctrain$, due to transport capacity issues. 
%(recall the problem setting in Section \ref{BlockCar}).
An artificial arc is defined for demand $k \in \cK$ from its last feasible \textsc{din} node, i.e., the last moment demand $k$ may leave its origin terminal and still get at its destination on time, and the container-sink node \textsc{csink} of its destination.

\section{The SSND-RM Formulation} \label{model}

The 
%\emph{Scheduled Service Network Design with Resource Management},
SSND-RM formulation for the intermodal rail blocking and railcar management problem includes the following design and flow decision variables:
%are defined in Table \ref{tab:decisions}. 
\begin{description}
    \item[$s_\sigma = 1$,] if train service $\sigma \in \Sigma_{\textsc{extra}}$ is selected, 0 otherwise;
    \item[$\yb = 1$,] if block $b \in \cB$ is selected, 0 otherwise;
    \item[$\zrk$,] number of containers of demand $k \in \cK$ on block $b \in \cBk$; 
    \item[$z_k$,] number of containers of demand $k \in \cK$ on the artificial arc associated with this demand; 
    \item[$x^\gamma_b$,] number of loaded railcars of type $\gamma \in \Gamma$ assigned to block $b$;
    \item[$w^\gamma_b$,] number of empty railcars of type $\gamma \in \Gamma$ assigned to block $b$;
    \item[$w^\gamma_\theta$,] number of empty railcars of type $\gamma \in \Gamma$ allocated to terminal $\theta \in \Theta$;
    \item[$w^\gamma_{\theta i}$,] number of empty railcars of type $\gamma \in \Gamma$ in pool node $i \in \cNpoolterm$ of terminal $\theta \in \Theta$; 
    \item[$\nu^{\tau, \tau^{\prime}}_{b \pi}$] number of platforms of type $\pi \in \Pi$ loaded with two containers of types $\tau$ and $\tau' \in \cT$\\ (no container type $\tau_{53}$ on bottom slot of $\pi_{40}$ platforms);
    \item[$\nu^{\tau}_{b \pi}$,] number of platforms of type $\pi \in \Pi$ loaded with a single container of type $\tau \in \cT$ \\ (no container type $\tau_{53}$ on $\pi_{40}$ platforms).
\end{description}

The ILP of our problem is:
\begin{equation} \label{obj1}
\mbox{Minimize} \sum_{\sigma \in \Sigma_{\textsc{extra}}} f_{\sigma} s_\sigma + \sum_{b \in \cB} \fb \yb + \sum_{b \in \cB}\sum_{k \in \cKr} c_{bk}z_{bk} + \sum_{k \in \cK} c_k \zk + \sum_{b \in \cB}\sum_{\gamma \in \Gamma} c_{b\gamma}w^\gamma_b + \sum_{\gamma \in \Gamma} \sum_{\theta \in \term} c_{\gamma\theta} w^\gamma_\theta 
\end{equation}
\noindent
Subject to:
%\emph{Demand block selection \& delivery}
\begin{flalign} 
\label{eq:model1-cstr-demand-left}
&\sum_{b \in \cBk} \zrk + \zk = \upsilon_k ,& &k \in \cK ,\\
\label{eq:model1-linking-zbk-yb}
&\zrk \leq \upsilon_k y_{b}  ,& &b \in \cBk , \ k \in \cK,\\
\label{eq:model1-cstr-loading-eq2}
    &\sum_{k \in \cK_b | \tk = \tau} \zrk = \sum_{\pi \in \Pi} \left(\nu^{\tau}_{b \pi} + 2 \nu^{\tau,\tau}_{b \pi} + \sum_{\tau'\in\cT|\tau' \neq \tau}\nu^{\tau,\tau'}_{b \pi}\right),& &\tau \in \cT, b \in \cB,\\
\label{eq:model1-cstr-loading-upperbound-eq1}
&\sum_{\gamma \in \Gamma} \eta^\gamma_{\pi} x^\gamma_{b} \geq  \sum_{\tau \in \cT} 
\left(  \nu^{\tau}_{b \pi} + \nu^{\tau,\tau}_{b \pi} + \frac{1}{2} \sum_{\tau^{\prime} \in \cT|\tau' \neq \tau} \nu^{\tau, \tau^{\prime}}_{b \pi} \right),& &\pi \in \Pi, b \in \cB ,\\
\label{eq:model1-cstr-loading-lowerbound-eq1}
&\sum_{\gamma \in \Gamma_\pi} x^\gamma_{b} \leq  \sum_{\tau \in \cT} 
\left(  \nu^{\tau}_{b \pi} + \nu^{\tau,\tau}_{b \pi} + \frac{1}{2} \sum_{\tau^{\prime} \in \cT|\tau' \neq \tau} \nu^{\tau, \tau^{\prime}}_{b \pi} \right),& &\pi \in \Pi, b \in \cB ,\\
\label{eq:model1-cstr-loading-eq3}
&\sum_{\gamma \in \Gamma} \left \lceil \frac{\eta^\gamma_{\pi_{40}}}{2} \right \rceil x^\gamma_b \geq \nu^{\tau_{40}, \tau_{53}}_{b \pi_{40}} ,& &b \in \cB,\\
%&{\color{blue}\sum_{\gamma \in \Gamma} \eta^\gamma_{\pi_{40}} x^\gamma_b \geq \nu^{\tau_{40} \tau_{53}}_{b \pi_{40}}} ,& &b \in \cB,\\
\label{eq:model1-linking-xb-wb}
&\sum_{\gamma \in \Gamma}\lambda_\gamma(x^\gamma_{b}+w^\gamma_b) \leq\ub \yb,& &b \in \cB,\\
%&\textcolor{blue}{\lambda_\gamma(x^\gamma_{b}+w^\gamma_b) \leq \min\{\lambda_\gamma H^\lambda,\ub\} \yb,& &b \in \cB,}\\
%\emph{Railcar flow conservation - Pool nodes}
\label{eq:model1-car-balance-pool}
&w^\gamma_{\theta i-1} +  \sum_{b \in \cBtermiin} ( w^\gamma_b + x^\gamma_b) = w^\gamma_{\theta i} + \sum_{b \in \cBtermiout} ( w^\gamma_b + x^\gamma_b) ,& &i \in \cN^{\textsc{pool}}_\theta , \gamma \in \Gamma , \theta\in \Theta,
\end{flalign}
\begin{flalign}
\label{eq:model1-car-allocation-first}
&w^\gamma_\theta =  w^\gamma_{\theta n} 
+ \sum_{b \in \cB_{\theta}^+|\alpha_b > \beta_b} ( w^\gamma_b + x^\gamma_b) ,& & \gamma \in \Gamma , \theta \in \Theta,\\
\label{eq:model1-car-allocation}
&\sum_{\theta \in \Theta} w^\gamma_\theta  \leq H^\gamma ,& &\gamma \in \Gamma ,\\
%linking cars, blocks, trains
\label{eq:model-cstr-capacity}
&\sum_{b \in \cBa} \sum_{\gamma \in \Gamma} \lambda_\gamma(x^\gamma_{b}+w^\gamma_b) \leq \ua,& &a \in \cAtrainm, \ \sigma \in \Sigma_\textsc{init},\\ 
\label{eq:model-cstr-extratrain-capacity}
&\sum_{b \in \cBa} \sum_{\gamma \in \Gamma} \lambda_\gamma(x^\gamma_{b}+w^\gamma_b) \leq s_\sigma \ua ,& &a \in \cAtrainm, \ \sigma \in \Sigma_\textsc{extra},\\
\label{eq:extraMinU}
&\frac{1}{u_a}\sum_{b \in \cBa} \sum_{\gamma\in\Gamma} \lambda_\gamma (x^\gamma_{b}+w^\gamma_b)
\geq U_\sigma s_\sigma,& &a \in \mathcal{A}_{\sigma}^{\textsc{utm}}, \sigma \in \Sigma_\textsc{extra}, \\ 
%\mathcal{A}_{\sigma}^{\text{minTM}}
%&\sum_{a \in \cAtrainm}\sum_{b \in \cBa} \sum_{\gamma\in\Gamma} d_a\lambda_\gamma (x^\gamma_{b}+w^\gamma_b) - \sum_{a \in \cAtrainm} (1+U) d_a \ua \nonumber & & \\ 
%&\geq - (1+U) \left(\max_{\sigma'\in\Sigma}\sum_{a \in \cA^{\textsc{tm}}_{\sigma'}} d_a \ua\right) p_\sigma,& &\sigma \in \Sigma_\textsc{extra},\\
%\label{eq:model-cstr-extratrain-bigM}
%&s_\sigma + p_\sigma \leq 1,& &\sigma \in \Sigma_{\textsc{extra}},\\
&y_b \in \{0,1\} ,& &b \in \cB, \nonumber\\
%	&y_{bk} \in \{0,1\} ,& &b \in \cB , \ k \in \cK,\\
&s_\sigma \in \{0,1\} ,& &\sigma \in \Sigma_{\textsc{extra}},\nonumber\\
%&p_\sigma \in \{0,1\} ,& &\sigma \in \Sigma_{\textsc{extra}},\nonumber\\
&x^\gamma_b, w^\gamma_\theta, w^\gamma_b, w^\gamma_{\theta i}, z_{bk}, \zk,  \nu^{\tau \tau^{\prime}}_{b\pi}, \nu^{\tau}_{b \pi} \in \mathbb{N},& &k \in \cK, b \in \cB , \tau \in \cT,\tau^{\prime} \in \cT, \pi\in \Pi,\nonumber\\
& & &i \in \cN, j \in \cN, \gamma \in\Gamma , \theta \in \Theta .\nonumber
\end{flalign}

The objective function (\ref{obj1}) minimizes the total cost of selecting extra trains (when relevant), selecting and building blocks, routing demand flows on blocks and artificial arcs, moving empty railcars, and allocating railcars to terminals initially.

Constraints~\eqref{eq:model1-cstr-demand-left} ensure that all the volume of each demand is transported either on blocks (train services) or on the artificial arcs. 
The linking constraints \eqref{eq:model1-linking-zbk-yb} enforce containers to be moved on selected blocks only, limiting the volume of any demand moving on a particular block to the total volume of that demand.
Constraints \eqref{eq:model1-cstr-loading-eq2} define the total container-flow (of a given type) moved by a particular block as the number of platforms (loaded with containers of that type) on the railcars within the block.
Notice that, jointly with the railcar-to-block linking constraints~\eqref{eq:model1-linking-xb-wb}, Constraints~\eqref{eq:model1-cstr-loading-eq2} enforce the flow of containers to be loaded on selected blocks only.

Recall that, a railcar is considered loaded as soon as at least one platform carries a container, even if the other platforms are empty.
Consequently, Constraints~\eqref{eq:model1-cstr-loading-upperbound-eq1} ensure that, the number of platforms, of a given type, making up the loaded railcars of a block is at least equal to the number of platforms of that type used on the block. 
Similarly, Constraints~\eqref{eq:model1-cstr-loading-lowerbound-eq1} ensure that the number of loaded railcars composed of platforms of a certain type on a block is not larger than the number of platforms of that type used on the block. 
Constraints~\eqref{eq:model1-cstr-loading-eq3} address the particular restrictions regarding the number of 40-on-bottom-53-on-top configurations on 40-ft platforms imposed by the physical rules of loading (see Section~\ref{BlockCar} for more details).
%\textcolor{red}{EF: I will make sure to explain this in the problem description. The max nb of plaforms with 53 on top is 2 on 3-platform railcars, and 3 on 5-platform railcars. Hence we need the divide by two and ceiling.}

Constraints~\eqref{eq:model1-linking-xb-wb} link the utilization of blocks by loaded and empty railcars to the selection of the block and its length.
%Let $\cBtermiout$ be the set of outgoing blocks made out of loaded (through arcs $\cApvl$ and $\cAcb$) and empty (arcs in $\cApeb$) railcars, at pool node $i \in \cNpoolterm$ of terminal $\theta \in \term$ (Figure \ref{fig:network}).
%Similarly, let $\cBtermiin$ be the sets of incoming blocks providing empty railcars to pool node $i \in \cNpoolterm$ of terminal $\theta \in \term$, through arcs $\cAbdp$ (empty in the block) as well as $\cAcd$ and $\cAcdp$ (loaded in the block and dismantled).
Constraints~\eqref{eq:model1-car-balance-pool} enforce the conservation of empty-railcar flows, of all types, at the pool nodes of the terminals.
Constraints~\eqref{eq:model1-car-allocation-first} count the number of railcars of each type needed at every terminal (initial allocation). The limits on the initial allocations are imposed through Constraints~\eqref{eq:model1-car-allocation}. Note that one may impose an initial allocation $H_{\theta}^{\gamma}$ to each terminal or let the tactical planning model decide on the best allocation, with or without a $H_{\theta}^{\gamma}$ limit.

Constraints~\eqref{eq:model-cstr-capacity} and \eqref{eq:model-cstr-extratrain-capacity} enforce the capacities of the regular and extra train services, respectively, on each of their legs. 
The latter also link the utilization of extra services to their selection.
Finally, constraints~\eqref{eq:extraMinU} enforce the condition that an extra train may be added to the schedule only if it is to be loaded at more than $U_\sigma\%$ of its capacity on a set of its legs. 
%\textcolor{red}{EF: This will be defined in the problem description.}
%{\color{magenta} \textbf{I defined the appropriate set of legs; I am assuming that this is what is in the code !!!}}

\section{Solution Approach} \label{solApp}

Before presenting our solution approach, we discuss key distinctive features of the SSND-RM formulation and outline the practical considerations that guided our work.

%{\color{magenta}NOTE Theo: ALL decision variables are \textbf{strongly} related through linking and such constraints; It is somewhat weak as justification.Besides, how does one define "complexity"for a decision variable? I am making a couple of suggestions in this and the partitioning paragraph below; text may be improved :-)}

The SSND-RM formulation is based on a time-space network whose discretization is governed by the train schedule. Except for the selection of extra trains, the schedule is fixed and given. The set of extra train selection variables is particularly challenging because it introduces an additional layer of design decisions that constrain block selection decisions.
The feasible set of blocks $\cB$ is large, but is restricted by the train schedule. In addition, the demand and railcar flow variables are integers and are strongly linked through the platform loading constraints. 

Given the tactical nature of the problem, an overnight computation time budget is considered acceptable. For deployment of a decision-support system based on the model, reliability and code maintenance are important considerations. For these reasons, we give priority to solution approaches using general-purpose MILP solvers as a black box, and algorithms with solution quality guarantees.

For our application, preliminary experiments using a general-purpose MILP solver based on an a priori enumeration of the set of block paths $\cB$ revealed large optimality gaps even after 12 hours of computation. However, it is well known that relatively simple heuristics combined with general-purpose MILP solvers can yield substantial computation speed-ups \citep[e.g.,][]{Fischetti2017,JoncourEtAl23}. Inspired by such ideas, we propose a construction heuristic close to the relax-and-fix procedure described in Section~3.6.1 in~\cite{PochetWolsey06}. We use the resulting feasible solution as a warm start for a general-purpose solver. Thus, the overall solution approach we propose remains exact.

We partition the decision variables into three disjoint sets. The first set contains the most challenging variables, in our case, the binary extra-train variables. The second set orders groups of decision variables in decreasing order of perceived impact on the structure of the solution,
$\mathcal{R} := (\{y_b, x_b^\gamma, b \in \cB, \gamma \in \Gamma\}, \{w_b^\gamma, w_\theta^\gamma, w_{\theta i}^\gamma; b \in \cB, \theta \in \Theta, i \in \cNpoolterm, \gamma \in \Gamma\})$. The third set groups container and platform assignment to block decision variables, $\mathcal{S} := \{z_{bk}, z_k, \nu^{\tau}_{b \pi}, \nu^{\tau, \tau^{\prime}}_{b \pi} ; b \in \cB, k \in \cK, \pi \in \Pi, \tau, \tau' \in \cT \}$.
%, as can be seen from Constraints~\eqref{eq:model1-cstr-loading-eq2}. 

% We partition decision variables into three disjoint sets. The first set, contains the most complicating variables, in our case the binary extra train variables. The second set orders sets of decision variables in decreasing order of complexity,  $\mathcal{R} := (\{y_b, x_b^\gamma, b \in \cB, \gamma \in \Gamma\}, \{w_b^\gamma, w_\theta^\gamma, w_{\theta i}^\gamma; b \in \cB, \theta \in \Theta, i \in \cNpoolterm, \gamma \in \Gamma\})$. The third set groups strongly linked decision variables, $\mathcal{S} := \{z_{bk}, z_k, \nu^{\tau}_{b \pi}, \nu^{\tau, \tau^{\prime}}_{b \pi} ; b \in \cB, k \in \cK, \pi \in \Pi, \tau, \tau' \in \cT \}$ and can be seen from Constraints~\eqref{eq:model1-cstr-loading-eq2}.
% We provide the details around the treatment of these sets in Algorithm~\ref{alg:relaxfix}. 

\begin{algorithm}[htbp]
    \caption{Compute warm-start solution}\label{alg:relaxfix}
    \begin{algorithmic}[1]
        \Require Extra train variables $\Sigma_\textsc{extra}$, ordered sequence of sets of decision variables $\mathcal{R}$, set 
        %set of strongly linked variables 
        $\mathcal{S}$, and precision $\varepsilon >0$, 
        %\State \text{Relax all variables}
        \State \emph{argrelax} $\gets$ \text{Solve continuous relaxation of SSND-RM} \Comment{Values of relaxed decision variables}
        \State \emph{restrictedModel} $\gets$ Fix to zero all variables in $\Sigma_\textsc{extra}$ whose value in \emph{argrelax} is less than $\varepsilon$; Impose integrality on the other variables in $\Sigma_\textsc{extra}$ %\Comment{Round and fix to zero}
        \State  \emph{argrestrict} $\gets$ Solve \emph{restrictedModel}
        \State \emph{restrictedModel} $\gets$ Fix all variables in $\Sigma_\textsc{extra}$ to their values in \emph{argrestrict} %\Comment{Round and fix}
        \ForAll {sets of decision variables $\mathcal{R}[i]$}
        \State \emph{restrictedModel} $\gets$ Fix to zero all variables in $\mathcal{R}[i]$ whose value in \emph{argrestrict} is less than $\varepsilon$; Impose integrality on the other variables in $\mathcal{R}[i]$ %\Comment{Round and fix to zero}
        \If{$i=last(\mathcal{R})$} 
        \State \emph{restrictedModel} $\gets$ Impose integrality on all variables in $\mathcal{S}$
        \EndIf
        \State  \emph{argrestrict} $\gets$ Solve \emph{restrictedModel}
        \EndFor
        \State \Return \emph{argrestrict}
    \end{algorithmic}
\end{algorithm}

Algorithm~\ref{alg:relaxfix} provides the details. The first five steps are dedicated to the extra train variables. 
Steps~6-10 concern the ordered sequence of decision variables $\mathcal{R}$ where
we apply an iterative round-and-fix procedure that only rounds down and fix the variables whose values are close to zero. 
In the last iteration, we impose integrality on the decision variables $\mathcal{S}$ for which the rounding-and-fixing scheme could cause feasibility issues (Steps~8-10). Finally, we note that the precision $\varepsilon$ should be close to zero; we use $1e^{-5}$.

\section{Numerical Results}
\label{application}

In this section, we report results using the large network of our industrial partner. Their network spans a large share of North America (see illustration in Figure~\ref{fig:CNnetwork} in the Appendix). To increase variability and protect confidentiality, we simulate instances based on a real historical train schedule and demand data. We describe the experimental setup in Section~\ref{experimentalSetup}, followed by an analysis of the results in Section~\ref{results}. We conclude by providing managerial insights based on the results.

\subsection{Experimental Setup}
\label{experimentalSetup}

%We base our experimental setup on the rail network of our industrial partner (as illustrated in Figure~\ref{fig:CNnetwork} in the Appendix, their network spans a large share of North America) and we generate problem instances using a real historical train schedule. 
In the following we provide details around our experimental setup.

\paragraph{Objective function.} The objective function is a specific case of~\eqref{obj1}. The associated parameter values have been carefully tuned in collaboration with subject matter experts of our industrial partner to ensure that the model generates realistic solutions. We refer to Appendix~\ref{app:ObjParams} for details. 

\paragraph{Extra Trains.} To evaluate the impact of extra trains we consider a setting where any of the existing train services in the schedule can be duplicated. Moreover, we fix the threshold $U_\sigma$ of minimum capacity to fill to select an extra train in~\eqref{eq:extraMinU} to 50\% for all $\sigma \in \Sigma_\textsc{extra}$. 

\paragraph{Block Generation.} We restrict the set of potential blocks $\cB$ to those that are considered feasible in practice. While generating the paths in a depth-first search, we impose practical constraints given by the train schedule. We note that for all instances it takes 51 seconds (the variance is close to zero) to generate all the blocks (around 47,000). Since this time is constant, it is not included in the computing times reported in Section~\ref{results}. 

\begin{figure}[htbp]
    \centering
    \subfloat{%
        \centering
    %\caption{Large1}
        \includegraphics[width=.49\linewidth]{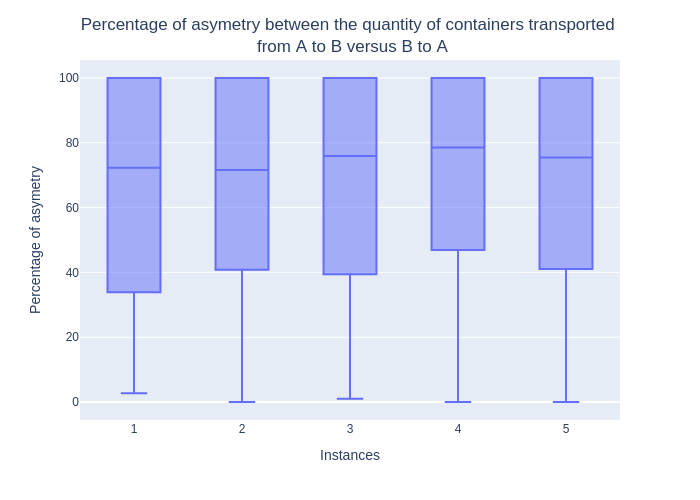}
        \label{fig:assymetry}}
    \subfloat{%
        \centering
    %\caption{Large2} 
        \includegraphics[width=.49\linewidth]{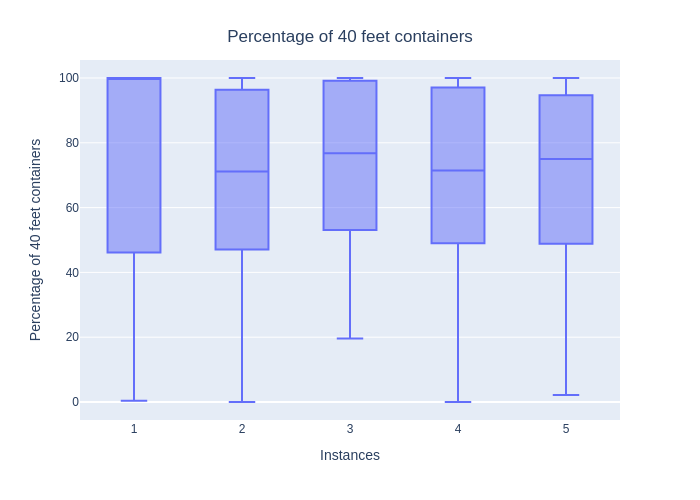}
        \label{fig:40ftContainers}}
    \caption{Illustration of how the demands are unbalanced and their types}
    \label{fig:Instance_illustration}
\end{figure}

\paragraph{Demand.} To keep real data confidential, we simulate five realistic demand instances based on actual values. There are 165 OD pairs in our network. To maintain the realism with respect to the scheduled train services (network capacity) that are fixed, we vary the characteristics of the demand while keeping overall demand volumes in the network intact. For this purpose, we randomly generate the proportion of 40-ft versus 53-ft containers for each OD pair. In Figure~\ref{fig:Instance_illustration}, we show that this leads to instances that have different characteristics. We display two box plots, each having one bar per instance (x-axis). The y-axis in the left-hand plot displays demand asymmetry. Based on total weekly volume for each pair of locations, it is computed as the percentage difference in volume (i.e., relative difference in demand between locations A to B, compared to B to A). 
The plot on the right-hand side, displays the percentage of 40-ft containers for each OD pair. First, we note that few OD pairs are perfectly symmetrical (zero percentage difference) and most are highly asymmetrical. Moreover, the percentage of 40-ft containers vary across ODs and across instances. 

\paragraph{Railcar fleet.} To assess the impact of different railcar types, we consider seven different railcar fleet scenarios:
\begin{enumerate}
    \item 1-platform 53-ft (1x53)
    \item 1-platform 40-ft, and 1-platform 53-ft (1x40,1x53)
    \item 3-platforms 53-ft (3x53)
    \item 3-platforms 40-ft, and 3-platforms 53-ft (3x40,3x53)
    \item 5-platforms 53-ft (5x53)
    \item 5-platforms 40-ft, and 5-platforms 53-ft (5x40,5x53)
    \item All six types (All)
\end{enumerate}
These scenarios cover different extreme cases: the most flexible cars in terms of loading but least efficient in terms of capacity usage (1x53), the most efficient in terms of capacity usage but most restricted in terms of loading (5x40), and the most flexible scenario from a network planning perspective (All).

\paragraph{Baseline Models.} We assess the value of the proposed model by comparing it to two baselines. First, a model that does not manage the railcar fleet hence assuming an infinite fleet, but that takes into account loading restrictions (called \emph{Unrestricted Fleet Model}). Second, a model that, in addition to ignoring the railcar fleet management, does not take into account loading constraints (called \emph{Unrestricted Loading Model}). The former model resembles the one proposed by \cite{Morganti2020}.
%, but it additionally incorporates the extra train service selection component. 
The latter model simply uses a formula that transforms a number of containers into capacity usage (train length) by using railcar lengths and assuming double stacking is always possible. We detail the integer linear formulations of the baseline models in Appendices~\ref{sec:unrestrictedFleet} and~\ref{sec:unrestrictedLoading}.

\paragraph{Solution Approaches and Hardware.} All experiments were performed using an Intel Core i9-10980XE 3.00 GHz processor with 128 GB of RAM. We used ILOG CPLEX 22.1.1, restricting the solver to a single thread, with a time limit of 12 hours and a stopping criteria of 2.5\% optimality gap. We compare the total solving times using CPLEX with our warm start approach to the time it takes to solve the SSND-RM formulation directly (without warm start).  

\paragraph{Performance metrics.} We use standard performance metrics: computing time in CPU seconds, optimality gap, root gap (optimality gap at the first found feasible solution) and number of instances that reach the time limit. 
%In addition, 
We provide interpretability by computing overall capacity utilization, percentage of railcar slot utilization, percentage of unsatisfied demand, fleet composition statistics, number of extra trains (if applicable), number of platforms, and number of railcars. These metrics are described in more detail next to the results. We solve five instances (demand scenarios) for each railcar scenario and version of the model/algorithm. Unless stated otherwise, we report an average over the five instances for each performance metric.

\subsection{Results}
\label{results}

The results are structured around the key takeaways. Namely, the warm-start solutions are of high quality, which leads to a significant speed up (Section~\ref{sec:computationalPerformance}), the baseline models are faster to solve but importantly underestimate the required capacity (Section~\ref{sec:baselines}), and railcar fleet management favours five-platform and 40-ft railcars whenever possible depending on the mix of container types to transport (Section~\ref{sec:res_railcarfleet}). 

%\ef{I have not yet updated the introduction to this section.} We first discuss our variable fixing iterative scheme in terms of optimality gaps and computing times. We follow suite by comparing those results (\textit{With Warmstart} approach) to the ones obtained by the \textit{Without Warmstart} one, further confirming the efficacy of our simple method. We then shift our focus on showing empirically the added value in taking the loading and fleet management problems into account by comparing three models: the one discussed in this paper (Fleet Management Model), the same one without fleet management (Unrestricted Fleet Model) and the same one without loading and fleet management (Unrestricted Loading Model). Finally, we go over the managerial insights that were found by analyzing the optimal solutions. 

\subsubsection{High-quality Warm-start Solutions Yield Significant Computational Speed-ups} \label{sec:computationalPerformance}

To assess computational performance, we solve the model with and without the extra train selection variables. 
%Moreover, we separate the results for the different railcar fleet scenarios. 
Table~\ref{table:general_results} reports the computing times for the different railcar fleet scenarios.
%comparing where each figure in the table is an average over five demand instances. WE SAID THAT 5 lines above !

Focusing on the left-hand side of the table, we note that it is relatively fast to compute the warm-start solutions. As expected, instances with extra trains are considerably more time consuming to solve than those without extra trains. The resulting warm-start solutions for the former instances are of high quality, which leads to a total average computing time within 3.25 hours. Note that the performance metrics consider instances that reached the time limit and we report the number of such instances in column (\#Timed-Out). Hence, the average optimality gap may be higher than the 2.5\% stopping criterion without reaching the 12 hour time limit on average (excluding such instances results in a total average time of 1,627
seconds for scenario 3x53, 998 seconds for scenario 3x40, 3x53, and 731 seconds for scenario 5x40, 5x53). Interestingly, 5x40, 5x53 (and not All) is the railcar scenario with the highest number of instances reaching the time limit. This may be due to symmetry issues.

In the case of no extra trains, the quality of the warm-start solutions is already below the 2.5\% optimality gap stopping criteria. We note that the averages are inflated by a few large values. The median time to compute a warm-start solution is 37 seconds for the model without extra trains, and 392 seconds for the model with extra trains.

Now, turning to the right-hand side of the table, we note that the warm start outperforms the general-purpose solver alone across all railcar scenarios, with and without extra trains. A total of 39 instances time-out, 20 in the extra train case. All of them concern a railcar scenario with a mix of 40-ft and 53-ft platform railcars.

\begin{table}[htbp]
\centering
\caption{Average computing times with and without warm start, for instances with and without extra trains}
\begin{adjustbox}{max width=\textwidth}
\begin{tabular}{l|rrrrr|rrrr}
\toprule
\multirow{2}{*}{\textbf{Scenario}} & \multicolumn{5}{c|}{\textbf{Extra Trains - With Warm start}} & \multicolumn{4}{c}{\textbf{Extra Trains - Without Warm start}} \\
\cmidrule(lr){2-6} \cmidrule(lr){7-10}
& \textbf{Warm start Time (s)} & \textbf{Total Time (s)} & \textbf{\#Timed-Out} & \textbf{Root Gap} & \textbf{Opt Gap} & \textbf{Time (s)} & \textbf{\#Timed-Out} & \textbf{Root Gap} & \textbf{Opt Gap} \\
\midrule
1x53            & 672   & 1,017  & 0  & 4.18\% & 2.40\% & 3,691  & 0 & 87.01\% & 2.42\% \\
1x40,1x53       & 3,435  & 13,313 & 0  & 3.44\% & 2.45\% & 43,257 & 5 & 97.36\% & 56.34\% \\
3x53            & 78    & 9,916  & 1 & 6.32\% & 2.49\% & 3,757 & 0  & 97.29\% & 2.42\% \\
3x40,3x53       & 4,707  & 17,820 & 2 & 3.24\% & 2.61\% & 43,257  & 5 & 97.51\% & 48.42\% \\
5x53            & 24    & 31  & 0   & 1.94\% & 1.94\% & 565  & 0  & 97.46\% & 1.94\% \\
5x40,5x53       & 368   & 26,125  & 3 & 3.28\% & 2.66\% & 43,257  & 5 & 97.57\% & 82.84\% \\
All             & 1,393  & 13,522 & 0 & 2.74\% & 2.38\% & 43,261 & 5 & 97.58\% & 97.58\% \\
\midrule
\multirow{2}{*}{\textbf{Scenario}} & \multicolumn{5}{c|}{\textbf{No Extra Trains - With Warm start}} & \multicolumn{4}{c}{\textbf{No Extra Trains - Without Warm start}} \\
\cmidrule(lr){2-6} \cmidrule(lr){7-10}
& \textbf{Warm start Time (s)} & \textbf{Total Time (s)} & \textbf{\#Timed-Out} & \textbf{Root Gap} & \textbf{Opt Gap} & \textbf{Time (s)} & \textbf{\#Timed-Out} & \textbf{Root Gap} & \textbf{Opt Gap} \\
\midrule
1x53            & 16   & 19  & 0  & 0.70\% & 0.70\% & 308  & 0  & 92.71\% & 1.32\% \\
1x40,1x53       & 38   & 68  & 0  & 2.26\% & 2.26\% & 3,691 & 4  & 95.16\% & 2.24\% \\
3x53            & 15   & 19  & 0  & 0.79\% & 0.79\% & 427  & 0  & 95.24\% & 0.82\% \\
3x40,3x53       & 1,317 & 1,479 & 0 & 2.07\% & 2.00\% & 43,257 & 5 & 97.18\% & 25.88\% \\
5x53            & 15   & 22  & 0  & 1.45\% & 1.44\% & 232  & 0  & 96.88\% & 1.97\% \\
5x40,5x53       & 159  & 184  & 0 & 1.50\% & 1.51\% & 43,257 & 5 & 97.56\% & 82.84\% \\
All             & 7,258 & 7,453 & 0 & 1.78\% & 1.78\% & 43,261 & 5 & 97.26\% & 97.24\% \\
\bottomrule
\end{tabular}
\end{adjustbox}
\label{table:general_results}
\end{table}

\subsubsection{Baseline Models Underestimate the Required Capacity} \label{sec:baselines}
As described in Section~\ref{experimentalSetup}, we assess the value of the SSND-RM formulation by comparing solutions with those obtained from two simpler baseline models.  We focus on assessing how they differ in the estimation of required capacity, hence we limit the comparison to instances without extra trains. 
%Moreover, 
None of those instances reached the time limit. The results are presented in Table~\ref{table:model_comparison} and we note that the computing time for SSND-RM is the same as the corresponding results in Table~\ref{table:general_results}. The rightmost sets of columns report results for the unrestricted fleet and unrestricted loading models, respectively. 

In addition to computing time, we report the percentage of unsatisfied demand (i.e., demand volumes assigned to artificial arcs compared to total demand volumes in the network), and the relative difference in network capacity usage for each model compared to the SSND-RM. We measure capacity in train length distance (train service capacity multiplied by the distance covered by each train leg). The capacity usage is then calculated as the train service used capacity multiplied by distance, divided by total network capacity in the same unit. The figures in the table correspond to the percentage difference compared to the capacity usage of the SSND-RM (which therefore has a usage difference of 0\%).

% Let $\lambda_a$ be the distance of an arc $a \in (\cAtrainm \cup \mathcal{A}_{\sigma}^{\textsc{utm}})$. The network capacity usage is computed as $\frac{100\sum_{\sigma \in \Sigma} \sum_{a \in (\cAtrainm \cup \mathcal{A}_{\sigma}^{\textsc{utm}})} \lambda_a \sum_{b \in \cBa} \sum_{\gamma\in\Gamma} \lambda_\gamma (x^\gamma_{b}+w^\gamma_b)}{\sum_{\sigma \in \Sigma} \sum_{a \in (\cAtrainm \cup \mathcal{A}_{\sigma}^{\textsc{utm}})} \lambda_a \ua}$. 
% The arc distances are taken into account as arcs with the same capacity should not have the same capacity weight if their length is different.

%\ef{Rewrite this part to focus on the implication of the results and how to interpret the numbers. Explain why it does not change in ULM.}
As expected, the SSND-RM model is harder to solve than the baseline models (its average computing times are considerably larger), but it accurately models the repositioning of empty railcars, as well as the loading restrictions specific to each railcar type. Ignoring the repositioning of empty railcars (as in the Unrestricted Fleet model), means that the train capacity required for this purpose is considered available to transport demand. This represents an average of 2-6\% of the network capacity. That is, using the unrestricted fleet model leads to an underestimation of the capacity needed by an average of 2-6\%, depending on the railcar scenario. 

Ignoring intermodal loading constraints, in addition to railcar fleet management, has a larger impact, as can be seen from the solutions to the Unrestricted Loading model. The average underestimation of capacity usage ranges from 17\% to 27\%, depending on the railcar scenario.

We note that the reduction in unsatisfied demand aligns with the underestimation of required capacity: the larger the underestimation, the more capacity appears available, enabling more demand to be fulfilled. Finally, we note that the usage difference in the \emph{Unrestricted Loading Model} varies as a result of changes in the capacity usage of the SSND-RM model, rather than fluctuations in its own usage. Indeed, as expected, it remains constant across railcar scenarios, as evidenced by the constant computing time and percentage of unsatisfied demand.

\begin{table}[htbp]
\centering
\caption{Comparison to model baselines on instances without extra trains (averages over five instances)}
\begin{adjustbox}{max width=\textwidth}
\begin{tabular}{l|rrr|rrr|rrr}
\toprule
\multirow{2}{*}{\textbf{Scenario}} & \multicolumn{3}{c|}{\textbf{SSND-RM (with warm start)}} & \multicolumn{3}{c|}{\textbf{Unrestricted Fleet Model}} & \multicolumn{3}{c}{\textbf{Unrestricted Loading Model}} \\
\cmidrule(lr){2-4} \cmidrule(lr){5-7} \cmidrule(lr){8-10}
& \textbf{Time (s)} & \textbf{Uns. Dem.} & \textbf{Usage Diff} 
                 & \textbf{Time (s)} & \textbf{Uns. Dem.} & \textbf{Usage Diff}
                 & \textbf{Time (s)} & \textbf{Uns. Dem.} & \textbf{Usage Diff} \\
\midrule
1x53         & 19.4 & 5.20\% & 0\%  & 3.7 & 5.12\% & 5\%  & 1.8 & 1.58\% & 27\% \\
1x40,1x53    & 68.3 & 1.64\% & 0\%  & 5.7 & 1.64\% & 6\%  & 1.8 & 1.58\% & 25\% \\
3x53         & 18.3 & 2.56\% & 0\%  & 4.1 & 2.55\% & 5\%  & 1.8 & 1.58\% & 25\% \\
3x40,3x53    & 1,478.6 & 0.63\% & 0\% & 27.1 & 0.62\% & 4\% & 1.8 & 1.58\% & 20\% \\
5x53         & 21.8 & 0.80\% & 0\%  & 4.1 & 0.80\% & 4\%  & 1.8 & 1.58\% & 20\% \\
5x40,5x53    & 184.4 & 0.55\% & 0\%  & 8.0 & 0.55\% & 4\%  & 1.8 & 1.58\% & 17\% \\
\text{All} & 7,453.4 & 0.56\% & 0\% & 11.9 & 0.56\% & 2\% & 1.8 & 1.58\% & 18\% \\
\bottomrule
\end{tabular}
\end{adjustbox}
\label{table:model_comparison}
\end{table}

%As such, we are able to see the value of computing our \textit{Fleet Management model}. Recall that we are devising a tactical plan of a length of a week that is repeated throughout the season. As such, calculation speed is not as crucial as for an operational plan. Waiting for a maximum of 2 hours (all railcar scenario) to emit such a plan is, thus, quite reasonable. What our \textit{Fleet Management model} gives in return is very valuable. First, it gives us a complete plan that accounts for the management of the railcar fleet. It allows us to know the optimal fleet composition for the season, where the railcars should initially be allocated in the network and how to reposition them so that the plan can be repeated over the horizon. Second, it allows us to have a plan that can be executed in real life because it gives us the real capacity that will be used in the network. As such, trains will not be overloaded ensuring that the plan can be followed as it is. As the maximum calculation time is of about 2 hours, using the \textit{Fleet Management model} is a no brainer.

\subsubsection{Railcar Fleet Management Favors Five-platform and 40-ft Railcars Whenever the Mix of Container Types Allows} \label{sec:res_railcarfleet}
Table~\ref{table:official_results} reports more detailed results on the selection of extra trains and the composition of the railcar fleet. More precisely, we report the average number of extra trains, the number of platforms and railcars, their slot utilization (the number of slots used compared to the total number of slots available on railcars), and the composition of railcar types in the selected fleet. Since the number of platforms varies across railcar scenarios, we separate the percentage of 53-ft platforms and the percentage of one-, three-, and five-platform railcars.

At our tactical level, the plan is repeated every week. Unsatisfied demand therefore accumulates unless extra resources are used. As expected, extra trains reduce the percentage of unsatisfied demand. Since we impose a minimum capacity usage for such trains, it nevertheless remains non-zero. We note that the number of extra trains varies depending on the railcar scenario.

Recall that railcar fleet management needs to strike the right balance between loading flexibility and capacity usage. The 53-ft one-platform railcars are the most flexible in terms of loading, while they make the least efficient use of capacity. On the other side of the spectrum, the 40-ft five-platform railcars are the most restrictive in terms of loading, but they make the most efficient use of capacity (measured by the number of containers per train foot). Relatively few extra trains are required when the model can choose among the most capacity-efficient railcars (two extra trains on average, compared to 14.2 on average when restricted to using 53-ft one-platform railcars).

Next, turning to the average slot utilization (\%slotUsed in Table~\ref{table:official_results}), we note that it is relatively stable across the different railcar scenarios. As expected, it decreases as the percentage of five-platform railcars increases in the fleet. When the model can select railcars of any type, on average, 70\% of the platforms are on five-platform cars (68\% in the case of extra trains). We note that whenever the model can select between 40-ft and 53-ft railcars, with only one exception, a majority of the platforms are 40-ft (the average percentage of 53-ft platforms varies between 24\% and 45\%, with one exception at 55\%). The exception concerns the extra train instances, because the demand satisfied by the extra trains is mostly 53-ft containers.

\begin{table}[htbp]
\centering
\caption{Results on railcar fleet composition (averages over five instances)}
\begin{adjustbox}{max width=\textwidth}
\begin{tabular}{l|rrrrrrrrrr}
\toprule
\multirow{2}{*}{\textbf{Scenario}} & \multicolumn{10}{c}{\textbf{Extra Trains}}\\
\cmidrule(lr){2-11}
& \textbf{Time (s)} & \textbf{Uns. Dem.} & \textbf{\#extra} & \textbf{nbPlatforms} & \textbf{nbRailcars} & \textbf{\%slotUsed} & \textbf{\%53Ftplat.} & \textbf{\%1plat.} & \textbf{\%3plat.} & \textbf{\%5plat.} \\
\midrule
1x53         & 1,017  & 0.43\% & 14.2 & 14,752 & 14,752 & 91\% & 100\% & 100\% & 0\% & 0\% \\
1x40,1x53    & 13,313 & 0.32\% & 5 & 14,726 & 14,726 & 90\% & 27\%  & 100\% & 0\% & 0\% \\
3x53         & 9,916  & 0.34\% & 9 & 15,043 & 5,014  & 84\% & 100\% & 0\%   & 100\% & 0\% \\
3x40,3x53    & 17,820 & 0.26\% & 2.2 & 14,860 & 4,953  & 83\% & 36\%  & 0\%   & 100\% & 0\% \\
5x53         & 31    & 0.27\% & 3 & 15,013 & 3,003  & 80\% & 100\% & 0\%   & 0\%   & 100\% \\
5x40,5x53    & 26,125 & 0.23\% & 2 & 15,275 & 3,055  & 78\% & 42\%  & 0\%   & 0\%   & 100\% \\
All          & 13,522 & 0.21\% & 2 & 15,080 & 3,993  & 79\% & 55\%  & 29\%  & 3\%   & 68\% \\
\midrule
\multirow{2}{*}{\textbf{Scenario}} & 
\multicolumn{10}{c}{\textbf{No Extra Trains}} \\
\cmidrule(lr){2-11}
& \textbf{Time (s)} & \textbf{Uns. Dem.} & \textbf{\#extra} & \textbf{nbPlatforms} & \textbf{nbRailcars} & \textbf{\%slotUsed} & \textbf{\%53Ftplat.} & \textbf{\%1plat.} & \textbf{\%3plat.} & \textbf{\%5plat.} \\
\midrule
1x53         & 19    & 5.20\% & 0 & 13,229 & 13,229 & 92\% & 100\% & 100\% & 0\%   & 0\% \\
1x40,1x53    & 68    & 1.64\% & 0 & 14,607 & 14,607 & 90\% & 24\%  & 100\% & 0\%   & 0\% \\
3x53         & 18    & 2.56\% & 0 & 13,955 & 4,652  & 85\% & 100\% & 0\%   & 100\% & 0\% \\
3x40,3x53    & 1,479  & 0.63\% & 0 & 14,561 & 4,854  & 83\% & 33\%  & 0\%   & 100\% & 0\% \\
5x53         & 22    & 0.80\% & 0 & 14,635 & 2,927  & 80\% & 100\% & 0\%   & 0\%   & 100\% \\
5x40,5x53    & 184   & 0.55\% & 0 & 14,931 & 2,998  & 79\% & 39\%  & 0\%   & 0\%   & 100\% \\
All          & 7,453  & 0.56\% & 0 & 14,758 & 3,664  & 81\% & 45\%  & 19\%  & 11\%  & 70\% \\
\bottomrule
\end{tabular}
\end{adjustbox}
\label{table:official_results}
\end{table}

\begin{figure}[H]
    \centering
    \subfloat{%
        \includegraphics[width=.49\linewidth]{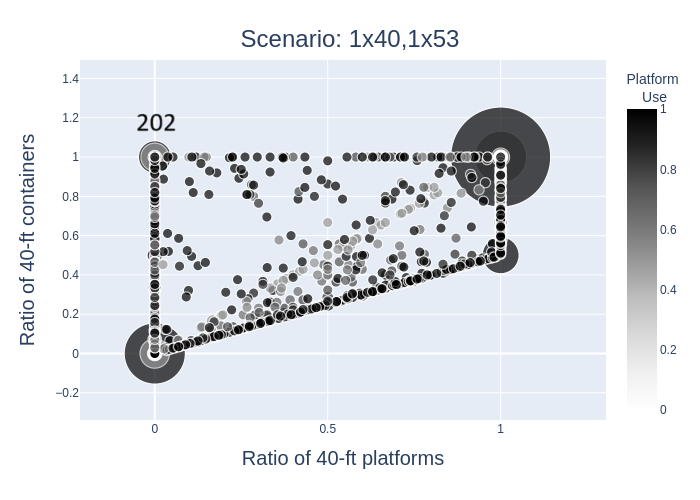}
        \label{fig:Scenario2}}
    \subfloat{%
        \includegraphics[width=.49\linewidth]{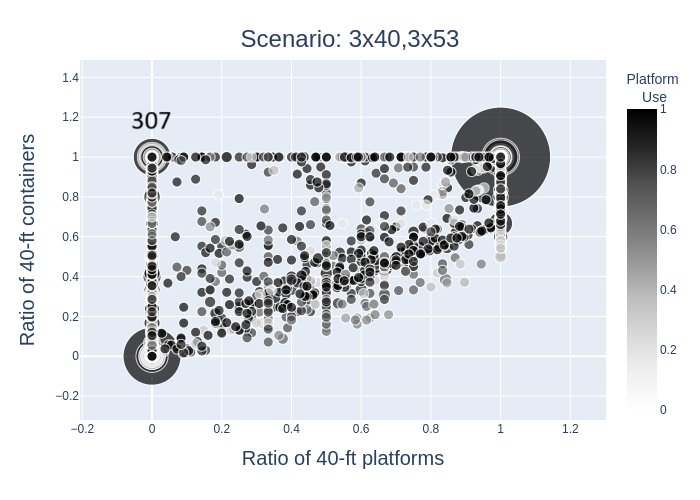}
        \label{fig:Scenario4}}
        \\
    \subfloat{%
        \includegraphics[width=.49\linewidth]{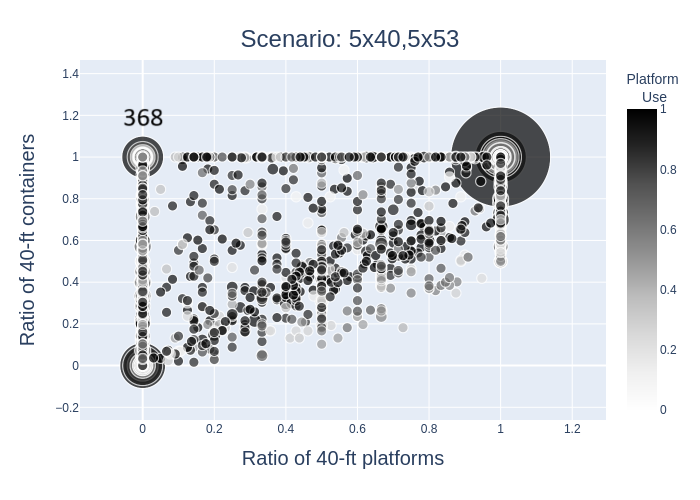}
        \label{fig:Scenario6}}
    \subfloat{%
        \includegraphics[width=.49\linewidth]{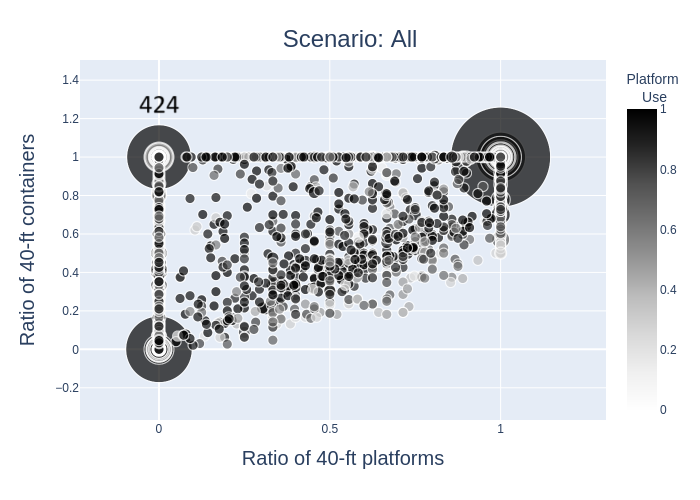}
        \label{fig:Scenario7}}
    \caption{Fleet composition for railcar scenarios with 40-ft and 53-ft platforms without extra trains}
    \label{fig:comp_illustration}
\end{figure}

Through the loading constraints, the composition of the railcar fleet is determined by the mix of container sizes in the demand volumes. To better analyze why, we provide a more granular view of the solutions in Figure~\ref{fig:comp_illustration}. To avoid instances that reached the time limit, we use only those without extra trains. The figure displays four scatter plots, one for each railcar scenario that allows selection between 40-ft and 53-ft railcars. Each dot represents a block in the five instances over which we averaged the statistics in the previous tables. The x-axis displays the ratio of 40-ft platforms in a given block, and the y-axis shows the ratio of 40-ft containers in the demand assigned to that block. The size and the shade of gray of each dot vary: the size increases with increasing density of dots in each plot. Dots of different shades may overlap, each dot’s size will vary according to the density of its respective shade. The legend for the shades of gray of the dot is shown on the right and represents the slot utilization, ranging from relatively low (light shade) to high (dark shade).

Each plot has a large dot at $(1,1)$, meaning that in each railcar scenario many blocks are made of 40-ft platforms that exclusively transport 40-ft containers. Similarly, each plot has a large dot at $(0,0)$, meaning that in each scenario many blocks are made of 53-ft platforms that transport exclusively 53-ft containers. The plots differ at $(0,1)$, which represents blocks made entirely of 53-ft platforms transporting exclusively 40-ft containers. This inefficient use of capacity is explained by unbalanced demand and the need to reposition railcars. The figure displayed above each of the dots at $(0,1)$ represents the total number of blocks, without taking into account the shade. Observe that, when the model is free to select any railcar type (bottom right), it sacrifices loading capacity efficiency (larger figure) in exchange for better use of the available platforms (larger black $(0,1)$ dot). Comparing the different plots at this point, we note that the total slot utilization improves when the model can select any railcar type.

The upper-left plot in Figure~\ref{fig:comp_illustration} displays a linear function between $(0,1)$ and $(1,0.5)$ that can be explained by the fact that the minimum ratio of 40-ft containers is half the ratio of the 40-ft platforms. Indeed, if only railcars with 40-ft platforms are available, at least half of the containers must be 40-ft, since 53-ft containers can only be loaded on top of 40-ft containers on 40-ft platforms. For all other railcar scenarios, some points are under this linear function, as some platforms can be empty.

%As shown by Table \ref{table:official_results} and \ref{fig:comp_illustration} , we can confidently say that the optimal fleet composition favors railcars from two specific types: railcars made of 5 platforms of 40 feet and 5 platforms of 53 feet.

%\begin{figure}[htb]
	%\centering
	%\includegraphics[width=0.5\textwidth]{Usage_large.png}
	%\caption{The given network provides a lot more capacity than needed}
	%\label{fig:usage}
%\end{figure}

\subsection{Managerial Insights}

Based on the results in the previous section, several general managerial insights emerge.

It is crucial to consider loading constraints for intermodal rail network planning. Ignoring such constraints can lead to a significant underestimation of the required capacity. In practice, this may lead to increased operational costs by requiring the addition of costly capacity on an ad hoc basis, and it may negatively impact customer service through delays.

The results also underline the importance of managing the railcar fleet. That is, accounting for the repositioning of empty railcars is important to adequately estimate network capacity usage. There is another side benefit: empty railcar repositioning represents unproductive moves of equipment. An accurate representation of such capacity can provide valuable information for adjusting the pricing strategy to better leverage this capacity and enhance profitability.

%The railcar fleet composition depends on demand characteristics (i.e., volumes and the mix of container types), and this link is captured through loading constraints. Block planning is less involved when ignoring railcar management and loading restrictions. In practice, such plans may be created by experienced humans. Our findings show that it is highly valuable to model loading restrictions and the repositioning of empty railcars, as the required network capacity can otherwise be underestimated. Since the plan is at a tactical level, underestimating the capacity may lead to an accumulation of demand, which may have implications for operations—for example, the necessity to add costly extra train services, which, in turn, can negatively impact network operations when not planned well in advance.

An adequate railcar fleet composition is important for the efficient use of network capacity. If demand in the network at a block level is characterized by a majority of 40-ft containers, then 40-ft railcars should be favored. For networks characterized by a high degree of variation in the mix of container types, the fleet should consist primarily of 53-ft platform railcars to ensure greater flexibility.

From a capacity usage point of view, there is a clear positive effect of using multi-platform railcars. Slot utilization per block is often used as a performance metric, especially to measure terminal productivity (loading at the origin of a block). Our results show that this metric alone may not be adequate. As an extreme example, a high average slot utilization can be achieved at a terminal level using one-platform 53-ft railcars, which may nevertheless result in poor overall network capacity utilization.

\section{Conclusion} \label{conclude}
We introduced a new problem for tactical planning of intermodal railway transport. The problem considers several interrelated decision-making problems capturing three consolidation processes and the management of a heterogeneous railcar fleet. We proposed a relatively simple solution approach that consists of generating warm-start solutions with a construction heuristic inspired by a relax-and-fix procedure. From a practical point of view, it has the advantage of leveraging a general-purpose solver as a black box.

Results based on realistic instances from our industrial partner, the Canadian National Railway Company, showed that we can solve large-scale instances down to an optimality gap of 2.5\% in a reasonable amount of time. We analyzed computational performance and the railcar fleet management solutions. Moreover, a comparison with solutions obtained by solving simpler baseline models showed that this new problem has high value. Based on the results, we drew managerial insights that are valuable to any railway company transporting intermodal containers on double-stack railcars.

This work opens up several avenues for future research. We solved a deterministic formulation of the problem, but at the tactical planning level there is typically uncertainty in the demand to be transported. A contextual stochastic formulation \citep{SADANA2025271} could lead to more resilient solutions. However, it poses significant computational challenges, especially since the demand distributions are expected to be decision-dependent \citep{Frejinger03072025}. Another important direction for future research is to tackle the corresponding operational network planning problem. In this case, there is less uncertainty about demand, but it is necessary to adapt the tactical plan to account for unforeseen events, such as equipment failures. In turn, this would require managing the railcar fleet at an even more granular level.

% Acknowledgments here
\section*{Acknowledgments}
We gratefully acknowledge the close collaboration with personnel from the Canadian National Railway Company, and the funding through the CN Chair in Optimization of Railway Operations at Universit\'e de Montr\'eal. Emma Frejinger was partially supported by the Canada Research Chair program [950-232244], T.G. Crainic by the NSERC Discovery Grant program, and Julie Kienzle was supported by a NSERC Graduate Scholarship [BESCD3-535243-2019]. While working on the project, T.G. Crainic held the UQAM Chair in Intelligent Logistics and Transportation Systems Planning and
was Adjunct Professor in the Department of Computer Science and Operations Research at the Universit\'e de Montr\'eal.
% Leave this (end of acknowledgment)

\bibliographystyle{abbrvnat}
\bibliography{railbib,IntermodalBlocking}

% Appendix here
% Options are (1) APPENDIX (with or without general title) or 
%             (2) APPENDICES (if it has more than one unrelated sections)
% Outcomment the appropriate case if necessary
%

\clearpage
\appendix 

\section{Intermodal network}

\begin{figure}[htbp]
    \centering
    \includegraphics[width=0.5\linewidth]{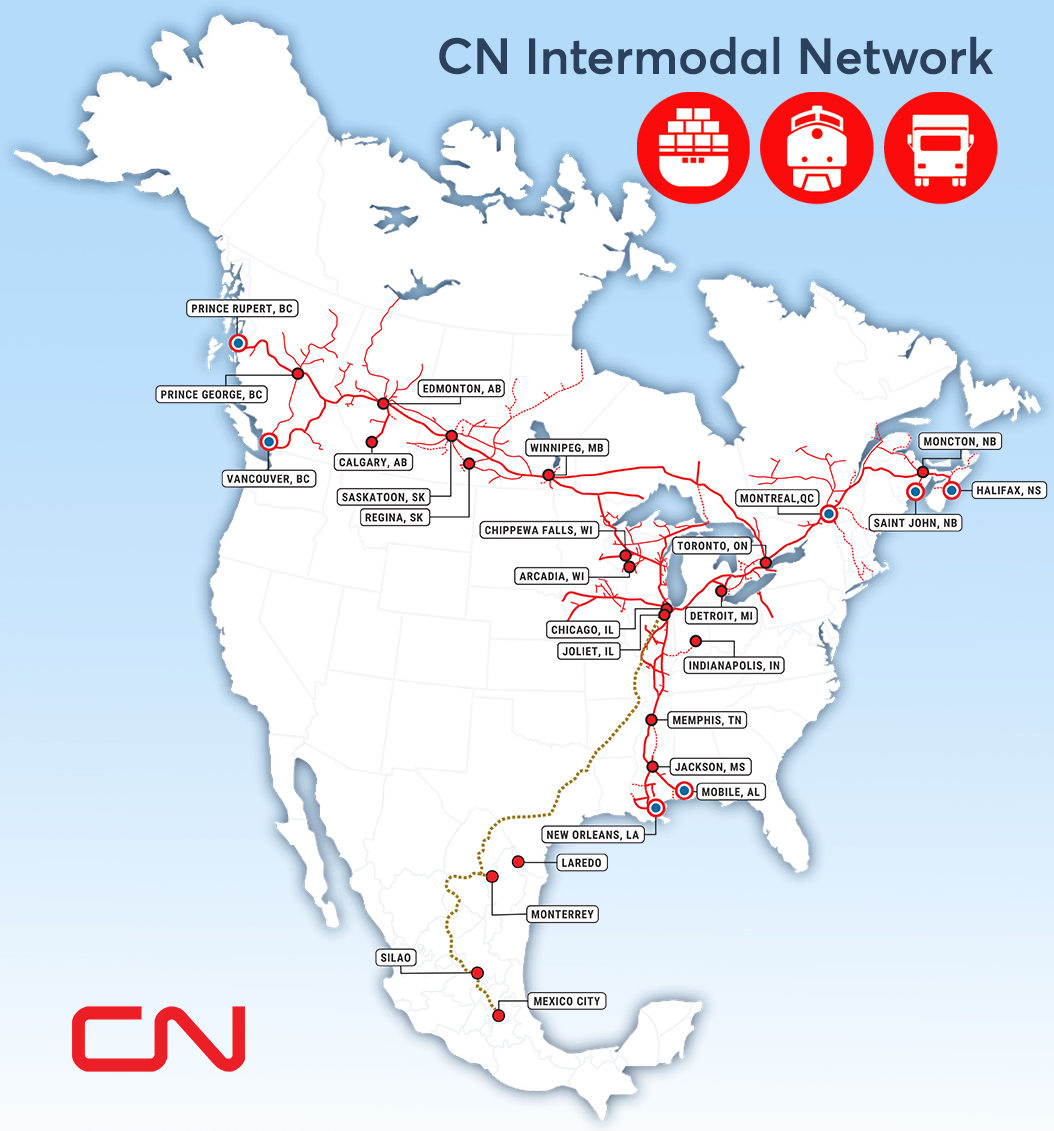}
    \caption{Canadian National Intermodal Network \citep{cn2025}}
    \label{fig:CNnetwork}
\end{figure}

\section{Objective Function Parameter Values}
\label{app:ObjParams}

%The model in the preceding sections is presented in a general way so that it can fit many different applications. In order to reach our objectives and run tests, some details specific to our application, such as the costs parameters in the objective function and the description of the set of extra trains that is used in our experiment, must be introduced.

%The model in the preceding sections was formulated generally to allow broad applicability. For the purposes of this study and the computational experiments, application-specific details are required, namely the cost parameters in the objective function and the set of additional trains.

The objective function for our application is the following:
\begin{align} \label{obj}
	&\mbox{Minimize} \sum_{b \in \cB} \left( \left(c^{\textsc{build}} + c^{\textsc{trans}} n^{\textsc{trans}}_b\right) \yb + \left(c^{\textsc{wait}}t^{\textsc{wait}}_b + c^{\textsc{bord}}n^{\textsc{bord}}_b + c^{\textsc{km}}d_b \right) \left(\sum_{k \in \cKr} z_{bk} + \sum_{\gamma \in \Gamma} w^\gamma_b\right)\right)\nonumber\\
	& + \sum_{b \in \cB}\sum_{k \in \cK} c^{\textsc{late}}t^{\textsc{late}}_{bk}z_{bk} + \sum_{\gamma \in \Gamma} \sum_{\theta \in \term} \eta^\gamma c^{\textsc{alloc}} w^\gamma_\theta + \sum_{k \in \cK} c^{\textsc{ndel}} \zk + \sum_{\sigma \in \Sigma_{\textsc{extra}}} \left(c^{\textsc{fix}} + c^{\textsc{var}}*\sum_{a\in \cAtrainm}\ua*d_a\right) s_\sigma
\end{align}

It incorporates several cost components reflecting operational considerations. A construction cost $c^{\textsc{build}}$ is incurred for each block built, while a transfer cost $c^{\textsc{trans}}$ is charged for every block transferred between trains ($n^{\textsc{trans}}_b$), penalizing the monetary and temporal burden of such operations. Transfers also generate a waiting cost $c^{\textsc{wait}}$ per empty railcar or demand for each minute of transfer time ($t^{\textsc{wait}}_b$). A cost $c^{\textsc{bord}}$ is incurred for every border crossing ($n^{\textsc{bord}}_b$). To discourage circuitous routes, a distance cost $c^{\textsc{km}}$ is applied per kilometer traveled 
($d_b$), complemented by a lateness cost $c^{\textsc{late}}$ per demand for each minute exceeding the shortest-path travel time ($t^{\textsc{late}}_{bk}$). An allocation cost $c^{\textsc{alloc}}$ per platform is introduced for assigning railcars at the start of the horizon, scaled by the number of platforms 
($\eta^\gamma$) of each railcar type $\gamma \in \Gamma$. To ensure high service quality, unmet demands (i.e., flows using artificial arcs) incur a large penalty cost $c^{\textsc{ndel}}$. Finally, additional trains incur a fixed cost $c^{\textsc{fix}}$ and a variable cost dependent on capacity and distance traveled. No explicit loading cost is included, as the loading process implicitly enforces maximum utilization of preferred platforms.

The cost values were established in close collaboration with our industrial partner.
%affiliated North American partner company. The parameters were selected to ensure the validity of the results, their consistency with the company’s operational practices and their suitability for computational implementation. 
They are the following: $c^{\textsc{build}} = 100$, $c^{\textsc{trans}} = 10020$, $c^{\textsc{wait}} = 1$, $c^{\textsc{bord}} = 1000$, $c^{\textsc{km}} = 0.75$, $c^{\textsc{late}} = 1$, $c^{\textsc{alloc}} = 200$, $c^{\textsc{ndel}} = 100 000$, $c^{\textsc{fix}} = 700 000$, $c^{\textsc{var}} = 0.01$. It is important to note that these do not represent the real monetary costs that the company pays to route demands through their railroad network.

\section{Unrestricted Fleet Model} \label{sec:unrestrictedFleet}

The model without fleet management can be written as:
\begin{flalign} \label{obj_model1}
\mbox{Minimize} \sum_{\sigma \in \Sigma_{\textsc{extra}}} f_{\sigma} s_\sigma + \sum_{b \in \cB} \fb \yb + \sum_{b \in \cB}\sum_{k \in \cKr} c_{bk}z_{bk} + \sum_{k \in \cK} c_k \zk&&
\end{flalign}
\noindent
Subject to:
%\emph{Demand block selection \& delivery}
\begin{flalign} 
\label{eq:model1-cstr-demand-left_model1}
&\sum_{b \in \cBk} \zrk + \zk = \upsilon_k ,& &k \in \cK ,\\
\label{eq:model1-linking-zbk-yb_model1}
&\zrk \leq \upsilon_k y_{b}  ,& &b \in \cBk , \ k \in \cK,\\
\label{eq:model1-cstr-loading-eq2_model1}
    &\sum_{k \in \cK_b | \tk = \tau} \zrk = \sum_{\pi \in \Pi} \left(\nu^{\tau}_{b \pi} + 2 \nu^{\tau,\tau}_{b \pi} + \sum_{\tau'\in\cT|\tau' \neq \tau}\nu^{\tau,\tau'}_{b \pi}\right),& &\tau \in \cT, b \in \cB,\\
\label{eq:model1-cstr-loading-upperbound-eq1_model1}
&\sum_{\gamma \in \Gamma} \eta^\gamma_{\pi} x^\gamma_{b} \geq  \sum_{\tau \in \cT} 
\left(  \nu^{\tau}_{b \pi} + \nu^{\tau,\tau}_{b \pi} + \frac{1}{2} \sum_{\tau^{\prime} \in \cT|\tau' \neq \tau} \nu^{\tau, \tau^{\prime}}_{b \pi} \right),& &\pi \in \Pi, b \in \cB ,\\
\label{eq:model1-cstr-loading-lowerbound-eq1_model1}
&\sum_{\gamma \in \Gamma_\pi} x^\gamma_{b} \leq  \sum_{\tau \in \cT} 
\left(  \nu^{\tau}_{b \pi} + \nu^{\tau,\tau}_{b \pi} + \frac{1}{2} \sum_{\tau^{\prime} \in \cT|\tau' \neq \tau} \nu^{\tau, \tau^{\prime}}_{b \pi} \right),& &\pi \in \Pi, b \in \cB ,\\
\label{eq:model1-cstr-loading-eq3_model1}
&\sum_{\gamma \in \Gamma} \left \lceil \frac{\eta^\gamma_{\pi_{40}}}{2} \right \rceil x^\gamma_b \geq \nu^{\tau_{40}, \tau_{53}}_{b \pi_{40}} ,& &b \in \cB,\\
%&{\color{blue}\sum_{\gamma \in \Gamma} \eta^\gamma_{\pi_{40}} x^\gamma_b \geq \nu^{\tau_{40} \tau_{53}}_{b \pi_{40}}} ,& &b \in \cB,\\
\label{eq:model1-linking-xb-wb_model1}
&\sum_{\gamma \in \Gamma}\lambda_\gamma x^\gamma_{b} \leq\ub \yb,& &b \in \cB,\\
%&\textcolor{blue}{\lambda_\gamma(x^\gamma_{b}+w^\gamma_b) \leq \min\{\lambda_\gamma H^\lambda,\ub\} \yb,& &b \in \cB,}\\
%\emph{Railcar flow conservation - Pool nodes}
%linking cars, blocks, trains
\label{eq:model-cstr-capacity_model1}
&\sum_{b \in \cBa} \sum_{\gamma \in \Gamma} \lambda_\gamma x^\gamma_{b} \leq \ua,& &a \in \cAtrainm, \ \sigma \in \Sigma_\textsc{init},\\ 
\label{eq:model-cstr-extratrain-capacity_model1}
&\sum_{b \in \cBa} \sum_{\gamma \in \Gamma} \lambda_\gamma x^\gamma_{b} \leq s_\sigma \ua ,& &a \in \cAtrainm, \ \sigma \in \Sigma_\textsc{extra},\\
\label{eq:model-cstr-extratrain-capacity_percentage_model1}
&\frac{1}{u_a}\sum_{b \in \cBa} \sum_{\gamma\in\Gamma} \lambda_\gamma x^\gamma_{b}
\geq U_\sigma s_\sigma,& &a \in \mathcal{A}_{\sigma}^{\textsc{utm}}, \sigma \in \Sigma_\textsc{extra},\\
%\mathcal{A}_{\sigma}^{\text{minTM}}
%&\sum_{a \in \cAtrainm}\sum_{b \in \cBa} \sum_{\gamma\in\Gamma} d_a\lambda_\gamma (x^\gamma_{b}+w^\gamma_b) - \sum_{a \in \cAtrainm} (1+U) d_a \ua \nonumber & & \\ 
%&\geq - (1+U) \left(\max_{\sigma'\in\Sigma}\sum_{a \in \cA^{\textsc{tm}}_{\sigma'}} d_a \ua\right) p_\sigma,& &\sigma \in \Sigma_\textsc{extra},\\
%\label{eq:model-cstr-extratrain-bigM}
%&s_\sigma + p_\sigma \leq 1,& &\sigma \in \Sigma_{\textsc{extra}},\\
&y_b \in \{0,1\} ,& &b \in \cB, \nonumber\\
%	&y_{bk} \in \{0,1\} ,& &b \in \cB , \ k \in \cK,\\
&s_\sigma \in \{0,1\} ,& &\sigma \in \Sigma_{\textsc{extra}},\nonumber\\
%&p_\sigma \in \{0,1\} ,& &\sigma \in \Sigma_{\textsc{extra}},\nonumber\\
&x^\gamma_b, z_{bk}, \zk,  \nu^{\tau \tau^{\prime}}_{b\pi}, \nu^{\tau}_{b \pi} \in \mathbb{N},& &k \in \cK, b \in \cB , \tau \in \cT,\tau^{\prime} \in \cT, \pi\in \Pi,\nonumber\\
& & &a \in \cA, i \in \cN, j \in \cN, \gamma \in\Gamma , \theta \in \Theta .\nonumber
\end{flalign}

The objective function (\ref{obj_model1}) minimizes the total cost of selecting extra trains (when relevant), selecting and building blocks, routing demand flows on blocks and artificial arcs.% and moving loaded railcars throughout the network.

Constraints~\eqref{eq:model1-cstr-demand-left_model1} to~\eqref{eq:model1-cstr-loading-eq3_model1} are the same as constraints~\eqref{eq:model1-cstr-demand-left} to~\eqref{eq:model1-cstr-loading-eq3} in the SSND-RM formulation. 
%Constraints~\eqref{eq:model1-car-balance-pool} to~\eqref{eq:model1-car-allocation} are removed as we do not account for the management of empty railcars. 
Constraints~\eqref{eq:model1-linking-xb-wb_model1} to~\eqref{eq:model-cstr-extratrain-capacity_percentage_model1} have the same role as Constraints~\eqref{eq:model1-linking-xb-wb}, and~\eqref{eq:model-cstr-capacity} to~\eqref{eq:extraMinU} in the SSND-RM formulation, but do not take into account the repositioning of empty railcars when calculating capacity usage.

\section{Unrestricted Loading Model} \label{sec:unrestrictedLoading}

The model without fleet management and loading is:
\begin{flalign} \label{obj_model2}
\mbox{Minimize} \sum_{\sigma \in \Sigma_{\textsc{extra}}} f_{\sigma} s_\sigma + \sum_{b \in \cB} \fb \yb + \sum_{b \in \cB}\sum_{k \in \cKr} c_{bk}z_{bk} + \sum_{k \in \cK} c_k \zk&&
\end{flalign}
\noindent
Subject to:
%\emph{Demand block selection \& delivery}
\begin{flalign} 
\label{eq:model1-cstr-demand-left_model2}
&\sum_{b \in \cBk} \zrk + \zk = \upsilon_k ,& &k \in \cK ,\\
\label{eq:model1-linking-zbk-yb_model2}
&\zrk \leq \upsilon_k y_{b}  ,& &b \in \cBk , \ k \in \cK,\\
%&{\color{blue}\sum_{\gamma \in \Gamma} \eta^\gamma_{\pi_{40}} x^\gamma_b \geq \nu^{\tau_{40} \tau_{53}}_{b \pi_{40}}} ,& &b \in \cB,\\
\label{eq:model1-linking-xb-wb_model2}
&\frac{1}{2}\sum_{\tau \in \cT}\sum_{k \in \cK_b | \tk = \tau}\lambda_\tau \zrk \leq \ub y_b,& &b \in \cB,\\
%&\textcolor{blue}{\lambda_\gamma(x^\gamma_{b}+w^\gamma_b) \leq \min\{\lambda_\gamma H^\lambda,\ub\} \yb,& &b \in \cB,}\\
%\emph{Railcar flow conservation - Pool nodes}
%linking cars, blocks, trains
\label{eq:model-cstr-capacity_model2}
&\frac{1}{2}\sum_{b \in \cBa} \sum_{\tau \in \cT}\sum_{k \in \cK_b | \tk = \tau}\lambda_\tau \zrk \leq \ua ,& &a \in \cAtrainm, \ \sigma \in \Sigma_\textsc{init},\\ 
\label{eq:model-cstr-extratrain-capacity_model2}
&\frac{1}{2}\sum_{b \in \cBa} \sum_{\tau \in \cT}\sum_{k \in \cK_b | \tk = \tau}\lambda_\tau \zrk \leq s_\sigma \ua ,& &a \in \cAtrainm, \ \sigma \in \Sigma_\textsc{extra},\\
\label{eq:model-cstr-extratrain-capacity_percentage_model2}
&\frac{1}{2u_a}\sum_{b \in \cBa} \sum_{\tau \in \cT}\sum_{k \in \cK_b | \tk = \tau}\lambda_\tau \zrk
\geq U_\sigma s_\sigma ,& &a \in \mathcal{A}_{\sigma}^{\textsc{tm}}, \sigma \in \Sigma_\textsc{extra},\\
%\mathcal{A}_{\sigma}^{\text{minTM}}
%&\sum_{a \in \cAtrainm}\sum_{b \in \cBa} \sum_{\gamma\in\Gamma} d_a\lambda_\gamma (x^\gamma_{b}+w^\gamma_b) - \sum_{a \in \cAtrainm} (1+U) d_a \ua \nonumber & & \\ 
%&\geq - (1+U) \left(\max_{\sigma'\in\Sigma}\sum_{a \in \cA^{\textsc{tm}}_{\sigma'}} d_a \ua\right) p_\sigma,& &\sigma \in \Sigma_\textsc{extra},\\
%\label{eq:model-cstr-extratrain-bigM}
%&s_\sigma + p_\sigma \leq 1,& &\sigma \in \Sigma_{\textsc{extra}},\\
&y_b \in \{0,1\} ,& &b \in \cB, \nonumber\\
%	&y_{bk} \in \{0,1\} ,& &b \in \cB , \ k \in \cK,\\
&s_\sigma \in \{0,1\} ,& &\sigma \in \Sigma_{\textsc{extra}},\nonumber\\
%&p_\sigma \in \{0,1\} ,& &\sigma \in \Sigma_{\textsc{extra}},\nonumber\\
& z_{bk}, \zk \in \mathbb{N},& &k \in \cK, b \in \cB , \tau \in \cT,\tau^{\prime} \in \cT, \pi\in \Pi,\nonumber\\
& & &a \in \cA, i \in \cN, j \in \cN, \gamma \in\Gamma , \theta \in \Theta .\nonumber
\end{flalign}

Objective function (\ref{obj_model2}) minimizes the total cost of selecting extra trains (when relevant), selecting and building blocks and routing demand flows on blocks and artificial arcs.

Constraints~\eqref{eq:model1-cstr-demand-left_model2} and~\eqref{eq:model1-linking-zbk-yb_model2} are the same as Constraints~\eqref{eq:model1-cstr-demand-left} and~\eqref{eq:model1-linking-zbk-yb} in the SSND-RM formulation. 
%Constraints~\eqref{eq:model1-cstr-loading-eq2} to~\eqref{eq:model1-cstr-loading-eq3} and~\eqref{eq:model1-car-balance-pool} to~\eqref{eq:model1-car-allocation} are removed as we do not account for the loading of containers onto railcars~\eqref{eq:model1-cstr-loading-eq2}-\eqref{eq:model1-cstr-loading-eq3} and the management of empty railcars~\eqref{eq:model1-car-balance-pool}-\eqref{eq:model1-car-allocation}. 
Constraints~\eqref{eq:model1-linking-xb-wb_model2} to~\eqref{eq:model-cstr-extratrain-capacity_percentage_model2} have the same role as constraints~\eqref{eq:model1-linking-xb-wb} and~\eqref{eq:model-cstr-capacity} to~\eqref{eq:extraMinU} in the SSND-RM formulation, but do not take into account the repositioning of empty railcars and loading of containers when calculating capacity usage. The formula transforms a number of containers into capacity usage (train length) by using railcar lengths and assuming double stacking is always possible.

\end{document}